\documentclass{gtpart}
\usepackage{amsmath}
\usepackage{paralist}
\usepackage{graphics} 
\usepackage{epsfig} 
\usepackage{graphicx}
\usepackage{epstopdf} 

\usepackage{tcolorbox}
\tcbuselibrary{skins}
\tcbuselibrary{breakable}
\usepackage{hyperref}
\definecolor{greenn}{rgb}{0,0.5,0}
\hypersetup{
    colorlinks=true,
    linkcolor=blue,
    filecolor=magenta,
    citecolor=greenn}

\usepackage[left]{lineno}

\usepackage[mathscr]{eucal} 
\usepackage{amsfonts}
\usepackage{sansmath}
\usepackage{float}
\usepackage{caption}
\usepackage{subcaption}
\usepackage{longtable}
\usepackage{tikz}
\usetikzlibrary{arrows}
\usepackage[section]{algorithm}
\usepackage{algpseudocode}
\usepackage{algorithmicx}

\usepackage{titlesec}
\usepackage{multirow}
\usepackage{enumitem}
\numberwithin{equation}{section}

\usepackage[marginratio=1:1,height=8.5in,width=6in,tmargin=1in]{geometry}

\usepackage{chngcntr}

\counterwithin{figure}{section}
\counterwithin{table}{section}

\newcommand{\RR}[1]{\mathbf{R}^{#1}}
\newcommand{\T}{T}
\newcommand{\dd}{\mathrm{d}}
\newcommand{\ii}{\mathbf{i}}
\newcommand{\OO}{\mathscr{O}}
\newcommand{\cSch}[1]{[\hspace{-0.065cm}[ #1 ]\hspace{-0.065cm}]}
\newcommand{\Cinf}[1]{\mathbf{\mathit{C}}^{\infty}_{#1}}
\newcommand{\bluecolor}[1]{{\color{blue} #1}}
\newcommand{\Input}[1]{\mbox{\textbf{Input:}\,\,}\parbox[t]{0.91\linewidth}{#1}}
\newcommand{\Output}[1]{\mbox{\textbf{Output:}\,\,}\parbox[t]{0.88\linewidth}{#1}}
\newcommand{\rreturn}[1]{\textbf{return}\,\,{#1}}
\newcommand{\CommentNew}[1]{\Comment{\mbox{\small\color{greenn}{#1}}}} 
\newcommand{\CommentCode}[1]{\quad\mbox{\small\color{greenn}{\# #1}}} 
\newcommand{\breakalg}{\algstore{bkbreak}\end{algorithmic}\end{algorithm}\begin{algorithm}[H]\vspace{0.1cm}\begin{algorithmic}[1]\algrestore{bkbreak}}
\newcommand{\ec}[1]{\mbox{$#1$}}

\newcommand{\LineCode}{\mbox{\tiny $>>>$}\,}

\providecommand{\bysame}{\makebox[3em]{\hrulefill}\thinspace}

\allowdisplaybreaks

\newtheorem{theorem}{Theorem}[section]

\newtheorem{lemma}[theorem]{Lemma}
\newtheorem{proposition}{Proposition}

\theoremstyle{definition}

\newtheorem{remark}[theorem]{Remark}
\newtheorem{example}[theorem]{Example} %

\title[On Computational Poisson Geometry II]{On Computational Poisson Geometry II: \\ Numerical Methods}

\author[M. Evangelista-Alvarado, J. C. Ru\'iz-Pantale\'on and P. Su\'arez-Serrato]{M. Evangelista-Alvarado, J. C. Ru\'iz-Pantale\'on and P. Su\'arez-Serrato}
\address{Instituto de Matem\'aticas, Universidad Nacional Aut\'onoma de M\'exico (UNAM), Mexico City, Mexico \newline \newline \href{https://www.aimsciences.org/article/doi/10.3934/jcd.2021012}{doi:10.3934/jcd.2021012} \hfill J. Comput. Dyn. \textbf{8}, 273--307 (2021)}

\begin{document}

\begin{abstract}
We present twelve numerical methods for evaluation of objects and concepts from Poisson geometry. We describe how each method works with examples, and explain how it is executed in code. These include methods that evaluate Hamiltonian and modular vector fields, compute the image under the coboundary and trace operators, the Lie bracket of differential 1--forms, gauge transformations, and normal forms of Lie--Poisson structures on \ec{\RR{3}}. The complexity of each of our methods is calculated, and we include experimental verifications on examples in dimensions two and three.
\end{abstract}

\maketitle

\tableofcontents




\bigskip



    \section{Introduction}

This is our second paper in a series on the computational aspects of Poisson geometry.
In the first paper we showed how fundamental concepts from Poisson geometry could be operationalized into symbolic code \cite{CompuPoisson}.
We also provided an associated Python module with our implementations, for ease of execution\footnote{Available via \url{https://github.com/appliedgeometry/poissongeometry}}.

A \emph{Poisson manifold} \cite{Poisson, WePG, Dufour, Kosmann, Camille} is a smooth manifold $M$ equipped with a contravariant skew--symmetric 2--tensor field $\Pi$, \emph{called Poisson bivector field}, satisfying the equation
    \begin{equation}\label{EcJacobiPi}
        \cSch{\Pi,\Pi} = 0,
    \end{equation}
with respect to the Schouten--Nijenhuis bracket \ec{\cSch{,}} for multivector fields \cite{Michor-08,Dufour}. Let \ec{m=\dim{M}}, and \ec{x = (U; x^{1}, \ldots, x^{m})} be local coordinates on $M$, then $\Pi$ has the following representation \cite{Lich-77,WeLocal}:
    \begin{equation}\label{EcPiCoord}
        \Pi = \tfrac{1}{2}\Pi^{ij}\frac{\partial}{\partial{x^{i}}} \wedge \frac{\partial}{\partial{x^{j}}}
            \ =   \sum_{1 \leq i < j \leq m} \Pi^{ij}\frac{\partial}{\partial{x^{i}}} \wedge \frac{\partial}{\partial{x^{j}}}
    \end{equation}
The functions \ec{\Pi^{ij}=\Pi^{ij}(x) \in \Cinf{U}} are called the coefficients of $\Pi$, and \ec{\{\partial/\partial{x^{i}}\}} is the canonical basis for vector fields on \ec{U \subseteq M}.

Poisson manifolds are generalizations of symplectic manifolds. A Poisson manifold can be thought of informally as a space that is foliated by symplectic leaves. It is then possible to define Hamiltonian dynamics relative to the symplectic forms defined on each leaf, via the Poisson bracket. Comprehensive treatments are available for interested readers \cite{Dufour, Camille}.

This geometric context provides a route that describes Hamiltonian dynamics rigorously. Applications of Hamiltonian dynamics are almost ubiquitous in every scientific domain. More specifically, uses of symplectic forms provide a formalism for diverse phenomena. This theory was introduced by Poisson himself to describe celestial mechanics \cite{Poisson}. It was then used by Dirac to investigate quantum mechanics \cite{Dirac}, and notably Kontsevich showed that Poisson manifolds admit deformation quantizations \cite{Kont}.

In this paper, we present twelve methods that allow for numerical computations of concepts from Poisson geometry. These are listed in the following Table \ref{table:Funs-Algos-Exes}, along with their respective algorithms, and a list of (not exhaustive) relevant references for each one. We have indicated with an asterisk (*) the six methods that depend explicitly on our module for symbolic computation \textsf{PoissonGeometry} \cite{CompuPoisson}.
    \begin{table}[H]
        \centering
        \caption{Our numerical methods, with their corresponding algorithms, and examples where they are used. The right column is an informal summary of the algorithmic complexities, computed and presented in detail in Section \ref{sec:ComplexityPerformance}.} \label{table:Funs-Algos-Exes}
        \resizebox{\textwidth}{!}{
     \begin{tabular}{|l|c|l|c|}
            \hline
        \multicolumn{1}{|c|} {\textbf{Method}} & \textbf{Algorithm} & \multicolumn{1}{c|}{\textbf{Examples}} & \textbf{Complexity} \\
            \hline
            \hline
        \hyperref[AlgNumBivector]{\textsf{num\_bivector\_field}}
        & \ref{AlgNumBivector}        & \cite{Dufour, Bullo, Kosmann, Camille} & O($m^2$) \\
            \hline
        \hyperref[AlgNumMatrixBivector]{\textsf{num\_bivector\_to\_matrix}}
        & \ref{AlgNumMatrixBivector}  & \cite{Dufour,Bullo,Kosmann,Camille} & O($m^2$) \\
            \hline
        \hyperref[AlgNumHamVF]{\textsf{num\_hamiltonian\_vf}}
        & \ref{AlgNumHamVF}           & \cite{Koz-95, Bullo, TV-19, Newton, HNN1} & O($m^2$) \\
            \hline
        \hyperref[AlgNumPoissonBracket]{\textsf{num\_poisson\_bracket}}
        & \ref{AlgNumPoissonBracket}  & \cite{JacobiMath, Dufour, Kosmann, Camille} & O($m^2$) \\
            \hline
        \hyperref[AlgNumSharp]{\textsf{num\_sharp\_morphism}}
        & \ref{AlgNumSharp}           & \cite{Dufour, Kosmann, Camille} & O($m^2$) \\
            \hline
        \hyperref[AlgNumCoboundary]{\textsf{num\_coboundary\_operator}}\,*
         & \ref{AlgNumCoboundary}     & \cite{Naka, Dufour, Kosmann, Poncin, MarcutSl2} & O($2^m$) \\
            \hline
        \hyperref[AlgNumModularVF]{\textsf{num\_modular\_vf}}\,*
        & \ref{AlgNumModularVF}       & \cite{Reeb2, Dufour, Kosmann, Miranda, Camille, MVallYu, Pedroza} & O($2^m$) \\
            \hline
        \hyperref[AlgNumCurl]{\textsf{num\_curl\_operator}}\,*
         & \ref{AlgNumCurl}           & \cite{GrabowskiFR, Damianou, Dufour, Camille} & O($2^m$) \\
            \hline
        \hyperref[AlgNumOneFormsB]{\textsf{num\_one\_forms\_bracket}}\,*
         & \ref{AlgNumOneFormsB}      & \cite{Dufour, Kosmann, Camille, Grabowski} & O($m^2$) \\
            \hline
        \hyperref[AlgNumGauge]{\textsf{num\_gauge\_transformation}}
         & \ref{AlgNumGauge}          & \cite{GaugeBursz, GaugeNaranjo, GaugeClass} & O($m^7$) \\
            \hline
        \hyperref[AlgNumNormal]{\textsf{num\_linear\_normal\_form\_R3}}\,*
        & \ref{AlgNumNormal}          & \cite{LiuXU-92, Naka, Ginzburg, Dufour, Sheng, Bullo, Camille, GaugeClass, MarcutSl2, Obook} & O($m$) \\
            \hline
        \hyperref[AlgNumFRatiu]{\textsf{num\_flaschka\_ratiu\_bivector}}\,*
        & \ref{AlgNumFRatiu}          & \cite{GrabowskiFR, Damianou, PabLef, PabloWrinFib, PabBott} & O($m^6$) \\
            \hline
            \hline
     \end{tabular}
        }
    \end{table}

The following diagram illustrates the internal functional dependencies of the methods available in  \textsf{NumPoissonGeometry} \footnote{Avalaible via \url{https://github.com/appliedgeometry/NumericalPoissonGeometry}}.
\begin{center}
    \begin{tikzpicture}[
      font=\rmfamily\footnotesize,
      every matrix/.style={ampersand replacement=\&, column sep=2cm, row sep=.23cm},
      source/.style={draw, thick, rounded corners, inner sep=.2cm},
      to/.style={->, >=stealth', shorten >=0.5pt, semithick},
      every node/.style={align=center}]

      \matrix{
        {}; \& \node[source] (gauge) {\hyperref[AlgNumGauge]{\textsf{num\_gauge\_transformation}}}; \\
        \node[source] (matrix) {\hyperref[AlgNumMatrixBivector]{\textsf{num\_bivector\_to\_matrix}}};
        \& \node[source] (hamilton) {\hyperref[AlgNumHamVF]{\textsf{num\_hamiltonian\_vf}}}; \\
        \node[source] (sharp) {\hyperref[AlgNumSharp]{\textsf{num\_sharp\_morphism}}};
        \& \node[source] (bracket) {\hyperref[AlgNumPoissonBracket]{\textsf{num\_poisson\_bracket}}}; \\
        \node[source] (formsbracket) {\hyperref[AlgNumOneFormsB]{\textsf{num\_one\_forms\_bracket}}}; \& {}; \\
        {} \& \node[source] (fratiu) {\hyperref[AlgNumFRatiu]{\textsf{num\_flaschka\_ratiu\_bivector}}}; \\
        \node[source] (numbivector) {\hyperref[AlgNumBivector]{\textsf{num\_bivector}}}; \& \node[source] (normal) {\hyperref[AlgNumNormal]{\textsf{num\_linear\_normal\_form\_R3}}}; \\
            };

      \draw[to] (matrix.east) -- (hamilton.west);
      \draw[to] (matrix) -- (sharp);
      \draw[to] (hamilton) -- (bracket);
      \draw[to] (sharp) -- (formsbracket);
      \draw[to] (matrix.east) --++(0:10mm)to[out=0, in=180] (gauge.west);
      \draw[to] (numbivector.east) --++(0:11mm)to[out=0, in=180] (fratiu.west);
      \draw[to] (numbivector.east) --++(0:9mm)to[out=0, in=180] (normal.west);
    \end{tikzpicture}\label{diagram}
\end{center}

The methods presented here have classical applications to Mechanics---see for example \cite{Bullo}---and also to recent advances in computer-aided techniques for determining normal forms for Hamiltonian systems \cite{CaraLoca20}. Furthermore, there has been a recent surge of interest in understanding Hamiltonian dynamics as uses of this theory start to appear in the data analysis and machine learning communities. Without attempting to be exhaustive, recent domains of application include; the development of Hamiltonian Monte Carlo techniques \cite{Duane}, applications of symplectic integration to optimization \cite{JordanPNAS}, inference of symbolic models from data \cite{Cranmer}, and the development of Hamiltonian Neural Networks \cite{HNN1, HNN2}.

Our work has been specifically designed to be compatible with popular machine learning frameworks, as our code can be integrated into NumPy, Pytorch, or TensorFlow environments. Moreover, as we rely on lattice meshes for our evaluations, our results will also be of interest to researchers that use finite element methods.

We hope to contribute an additional dimension to the understanding of Poisson geometry, enabling everyone to carry out numerical experiments with our freely available open-source code, motivating the expansion of techniques that have so far been incorporated. To the best of our knowledge this is the first comprehensive implementation of these methods.

Our numerical techniques can inform and complement researchers' intuition and provide further insights. For example, if a certain vector field is not trivial it can imply that a given Poisson structure is not unimodular, and this can be verified numerically. Moreover, our module only needs the algebraic expression for the Poisson bivector to carry out this verification (and not the complete, explicit vector field).

We have included specific examples for each of the methods above where numerical computations would be desirable, or are relevant in published work. We also strongly believe that our algorithms in this paper can be useful in the following related fields.

Our methods could aid in the development of deep learning systems for flows on tori and spheres \cite{Cranmer2020}, in particular for the groups \ec{U(1) \cong \mathbf{S}^{1}} and  \ec{SU(2) \cong \mathbf{S}^{3}}. In a related direction, normalized flows on Lie groups have been investigated with respect to distributions of data points \cite{Falorci2019}.

We expect our methods to be useful for evaluating inference problems related to the $n$--body problem, which has been approached recently \cite{Battaglia2016}.

Hamiltonian dynamics have been used to learn systems of simulations \cite{Battaglia2019}. This innovative work pioneers deep learning techniques that combine Hamiltonian dynamics and integrators of Ordinary Differential Equations. With our methods, such an overarching program could be extended to include more diverse equations, such as the ones that define a Poisson structure, or that verify conservation of specific quantities, for instance, a unimodular flow.

Normalized flows and techniques from symplectic geometry have been used to find canonical Darboux coordinates \cite{Li2020}. This highlights that the use of symplectic, and more generally, Poisson techniques is computationally efficient because these structures preserve volume.

Models that learn to respect conservative laws, through the use of Hamiltonian dynamics have been recently developed \cite{HNN1}. These could also be extended to more ample classes using our methods.

Hamiltonian functions themselves, relative to a symplectic form, have been modelled using recurrent neural networks \cite{chen2019}. We expect extensions of these works to be possible using our methods as well, where now the functions may be Hamiltonian relative to a Poisson structure. Moreover, systems that learn the dynamics from observed orbit trajectory data could be potentially designed.

There are two more applications that are relevant to the computation of statistical quantities on Lie groups. First, Hamiltonian flows that are equivariant with respect to an action can be used to learn densities that are invariant under the action of a Lie algebra. This strategy has been realized already for very specific algebras \cite{Jimenez2019}, and could be extended to more with our methods. Second, these ideas could also lead to a general method for finding reparametrizable densities on arbitrary Lie groups, using their associated Lie-Poisson structure \cite{Tamayo2015}.

\medskip

This paper is organized in the following way. In section 2, we describe our algorithms and point to domains of application for each one. In section 3,  detailed analyses of the complexities of our methods are presented, as well as experimental run-times. A summary with explicit upper bounds is found in Table \ref{table:complexity}.

    \section{Numerical Methods} \label{sec:keyf}

|In this section we describe the implementation of all functions in \textsf{NumPoissonGeometry} and present numerical examples with classical/relevant Poisson bivector fields. Throughout this paper we will assume that the input functions used in our methods are computable.

    \subsection{\textsf{NumPoissonGeometry}: Evaluation} \label{subsec:mesh}

Our methods work with regular and \emph{irregular} meshes that must be entered as \textsf{NumPy} arrays: a \ec{(k,m)} \textsf{NumPy} array for a mesh with $k$ points in \ec{\RR{m}}, for \ec{k \in \bold{N}}.

Irregular meshes can be used to implement probabilistic/statistical methods in Poisson geometry. In several of the subsections below we generate irregular meshes by means of random samples drawn from  a uniform distribution over the interval \ec{[0,1)}. Sometimes for simplicity we use the `corners' of the unit cube in \ec{\RR{m}},
    \begin{equation}\label{EcCorners}
        Q^{m} :=  \underset{m\text{\ times}}{\underbrace{\{0,1\}\ \times \cdots \times \{0,1\}}},
    \end{equation}
preloaded in Python into a \textsf{NumPy} array called \textsf{Qmesh}.

    \subsection{\textsf{NumPoissonGeometry}: Syntax} \label{subsec:class}

The instance and syntax of \textsf{NumPoissonGeometry} are the same as that of \textsf{PoissonGeometry}. We recall these here for the reader's convenience. \\

\noindent\textbf{Coordinates.} By default, to emulate canonical coordinates on \ec{\RR{m}}, we use symbolic variables that are just the juxtaposition of the symbol \textsf{x} and an index of the set \ec{\{1,\ldots,m\}}: (\text{\textsf{x1}}, \ldots, \text{\textsf{xm}}). \\

\noindent\textbf{Scalar Functions.} A scalar function is written using \emph{string literal expressions}. For example, the function \ec{f = ax_1^2 + bx_2^2 + cx_3^2} should be written exactly as follows: \textsf{`a*x1**2 + b*x2**2 + c*x3**2'}. All characters that are not coordinates are treated as (symbolic) parameters: \textsf{a}, \textsf{b} and \textsf{c} for the previous example. \\

\noindent\textbf{Multivector Fields and Differential forms.} Both multivector fields and differential forms are written using \emph{dictionaries} with \emph{tuples of integers} as \textsl{keys} and \emph{string} type \textsl{values}. If the coordinate expression of an $a$--multivector field $A$, with \ec{a \in \mathbb{N}}, is given by
    \begin{equation*}
        A = \sum_{1 \leq i_1 < i_2 < \cdots < i_a \leq m} A^{i_1 i_2 \cdots i_a}\frac{\partial}{\partial{x_{i_1}}} \wedge \frac{\partial}{\partial{x_{i_2}}} \wedge \cdots \wedge \frac{\partial}{\partial{x_{i_a}}}, \quad A^{i_1 \cdots i_a} = A^{i_1 \cdots i_a}(x),
    \end{equation*}
then $A$ should be written using a dictionary, as follows:
    \begin{equation}\label{EcMultivectorDic}
        \Big\{(1,...,a): \mathscr{A}^{1 \cdots a}, \,...,\, (i_1,...,i_a): \mathscr{A}^{i_1 \cdots i_a}, \,...,\, (m-a+1,...,m): \mathscr{A}^{m-a+1\cdots m}\Big\}
    \end{equation}
Here, each key \ec{(i_1,\ldots,i_a)} is a tuple containing ordered indices \ec{1 \leq i_1 < \cdots < i_a \leq m} and the corresponding value \ec{\mathscr{A}^{i_1 \cdots i_a}} is the string expression of the scalar function (coefficient) \ec{A^{i_1 \cdots i_a}} of $A$.

The syntax for differential forms is the same. It is important to remark that we can only write the keys and values of \emph{non--zero coefficients}. See the documentation for more details.

    \subsection{\textsf{NumPoissonGeometry}: Python Implementation} \label{sec:implementation}

First let us briefly describe the syntax of the module \textsf{NumPoissonGeometry}. The inputs of the twelve methods in Table \ref{table:Funs-Algos-Exes} have to be string literal expressions (or dictionaries containing them). The \textsf{sympify} method converts such string expressions into symbolic variables, and the \textsf{lambdify} method transforms the symbolic expressions into functions that allow a (fast) numerical evaluation.

The output of our methods can be chosen to be a \textsf{NumPy} array (by default) or a \textsf{pytorch}/\textsf{tensorflow} tensor containing the evaluation of the input data at each point in a mesh on \ec{\RR{m}}.

Next we will explain each of our numerical methods, how they relate to the theoretical objects being implemented, present the corresponding algorithm, and show how they may be used with helpful examples.

    \subsection{Bivector Fields} \label{subsec:bivector}

The evaluation of a (Poisson) bivector field $\Pi$ at a point in $M$ is defined by
    \begin{equation*}
        \Pi_{p} := \tfrac{1}{2}\Pi^{ij}(p)\left.\frac{\partial}{\partial{x^{i}}}\right|_{p} \wedge \left.\frac{\partial}{\partial{x^{j}}}\right|_{p}, \quad p \in M.
    \end{equation*}
Observe that the coefficients of $\Pi$ in (\ref{EcPiCoord}) at $p$ determines the evaluation above.

The function \textsf{num\_bivector} evaluates a (Poisson) bivector field on a mesh in \ec{\RR{m}}.

\begin{algorithm}[H]
    \captionsetup{justification=centering}
    \caption{\ \textsf{num\_bivector}(\emph{bivector, mesh})} \label{AlgNumBivector}
        \rule{\textwidth}{0.4pt}
    \Input{a (Poisson) bivector field and a mesh}
    \Output{evaluation of the bivector field at each point of the mesh}
        \rule{\textwidth}{0.4pt}
    \begin{algorithmic}[1] 
        \Procedure{}{}
            \State $m$ $\gets$ dimension of the manifold
            \State \bluecolor{bivector} $\gets$ a variable encoding the (Poisson) bivector field
            \State \bluecolor{mesh} $\gets$ a $(k,m)$ array encoding the mesh
                \CommentNew{$k$: number of points in the mesh}
            \State \textsc{Transform} each \bluecolor{bivector} item into a function that allows a numerical evaluation
            \State \rreturn{an array containing the evaluation of \bluecolor{bivector} at each $m$--array of \bluecolor{mesh}}
        \EndProcedure
    \end{algorithmic}
\end{algorithm}

\begin{example}\label{example:so3}
Consider the Lie--Poisson bivector field on \ec{\RR{3}_{x}}
    \begin{equation}\label{EcPiSO3}
        \Pi_{\mathfrak{so}(3)} =
        x_3\frac{\partial}{\partial x_1}\wedge \frac{\partial}{\partial x_2} -
        x_2\frac{\partial}{\partial x_1}\wedge \frac{\partial}{\partial x_3} +
        x_1 \frac{\partial}{\partial x_2}\wedge \frac{\partial}{\partial x_3},
    \end{equation}
associated to the 3--dimensional Lie algebra \ec{\mathfrak{so}(3)}, which is used in studies of the 3--body problem \cite{Newton}, and more generally in geometric mechanics \cite{Bullo}. To evaluate \ec{\Pi_{\mathfrak{so}(3)}} at points of \ec{Q^{3}} in (\ref{EcCorners}) we compute:

\begin{tcolorbox}[arc=0mm, boxsep=0mm, skin=bicolor, colback=pink!15, colframe=blue!20, colbacklower=blue!0, breakable, halign=left]
\LineCode \textsf{npg3 = NumPoissonGeometry(3)}
    \hspace*{\fill} \CommentCode{\textsf{NumPoissonGeometry} instance} \\
\LineCode \textsf{P\_so3 = \big\{(1, 2): `x3', (1, 3): `-x2', (2, 3): `x1'\big\}}
    \hspace*{\fill} \CommentCode{dictionary for \ec{\Pi_{\mathfrak{so}(3)}} in (\ref{EcPiSO3}) according to (\ref{EcMultivectorDic})} \\
\LineCode \textsf{npg3.num\_bivector(P\_so3, Qmesh, dict\_output=True)}
    \hspace*{\fill} \CommentCode{run \textsf{num\_bivector} function}
    \tcblower
\textsf{[\parbox[t]{\linewidth}{\{(1, 2): 0.0, (1, 3): \phantom{-}0.0, (2, 3): 0.0\}, \
 \{(1, 2): 1.0, (1, 3): \phantom{-}0.0, (2, 3): 0.0\}, \newline
 \{(1, 2): 0.0, (1, 3): -1.0, (2, 3): 0.0\}, \
 \{(1, 2): 1.0, (1, 3): -1.0, (2, 3): 0.0\}, \newline
 \{(1, 2): 0.0, (1, 3): \phantom{-}0.0, (2, 3): 1.0\}, \
 \{(1, 2): 1.0, (1, 3): \phantom{-}0.0, (2, 3): 1.0\}, \newline
 \{(1, 2): 0.0, (1, 3): -1.0, (2, 3): 1.0\}, \
 \{(1, 2): 1.0, (1, 3): -1.0, (2, 3): 1.0\}]}}
\end{tcolorbox}

Note that the output preserves the \textsf{PoissonGeometry} syntax (\ref{EcMultivectorDic}). For example, to produce a \textsf{PyTorch} tensor encoding this information we use the \textsf{pt\_output} flag:
\begin{tcolorbox}[arc=0mm, boxsep=0mm, skin=bicolor, colback=pink!15, colframe=blue!20, colbacklower=blue!0, breakable, halign=left]
\LineCode \textsf{npg3.num\_bivector(P\_so3, Qmesh, pt\_output=True)}
    \hspace*{\fill} \CommentCode{run \textsf{num\_bivector} function with \textsf{pt\_output} flag}
    \tcblower
\resizebox{\textwidth}{!}{
\textsf{tensor([\,\parbox[t]{\linewidth}{[[0.,  0., \phantom{-}0.], [0.,  0.,  0.], [0., \phantom{-}0.,  0.]], \
        [[0.,  1., \phantom{-}0.], [-1.,  0.,  0.], [0., \phantom{-}0.,  0.]], \newline
        [[0.,  0., -1.], [0.,  0.,  0.], [1., \phantom{-}0.,  0.]], \
        [[0.,  1., -1.], [-1.,  0.,  0.], [1., \phantom{-}0.,  0.]], \newline
        [[0.,  0., \phantom{-}0.], [0.,  0.,  1.], [0., -1.,  0.]], \
        [[0.,  1., \phantom{-}0.], [-1.,  0.,  1.], [0., -1.,  0.]], \newline
        [[0.,  0., -1.], [0.,  0.,  1.], [1., -1.,  0.]], \
        [[0.,  1., -1.], [-1.,  0.,  1.], [1., -1.,  0.]]\,], \newline
        dtype=torch.float64)}}
}
\end{tcolorbox}
\end{example}

\begin{figure}[H]
\centering
    \begin{subfigure}{0.35\textwidth}
        \includegraphics[width=\textwidth]{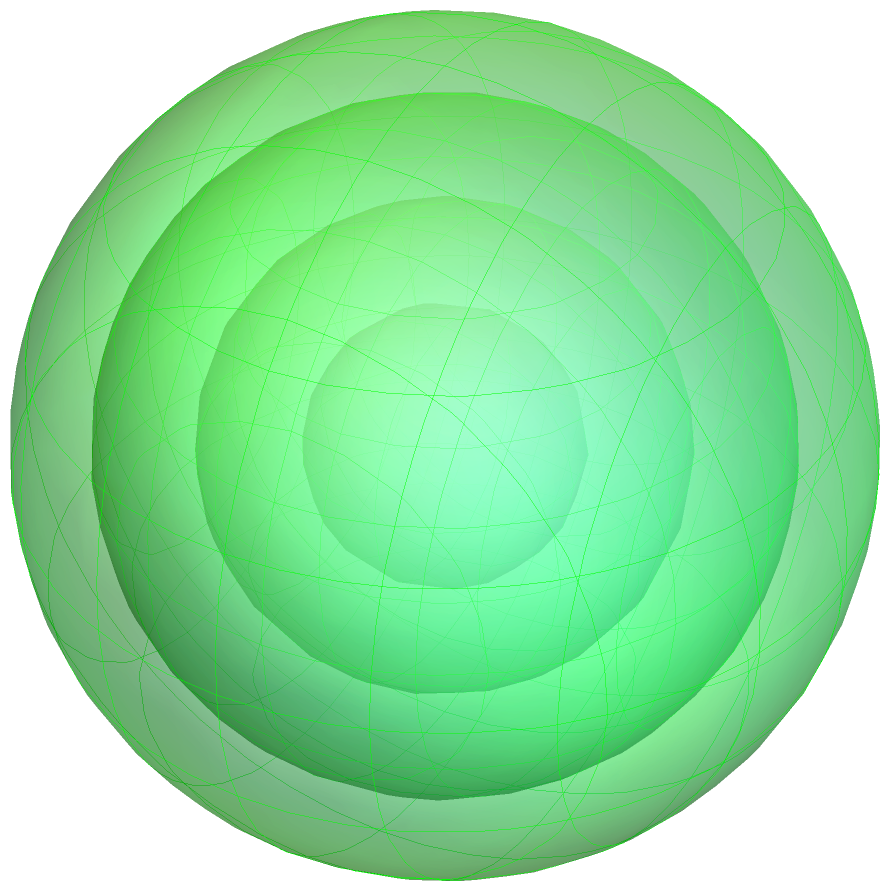}
    \end{subfigure}
        \hspace{0.6cm}
    \begin{subfigure}{0.4\textwidth}
        \includegraphics[width=\textwidth]{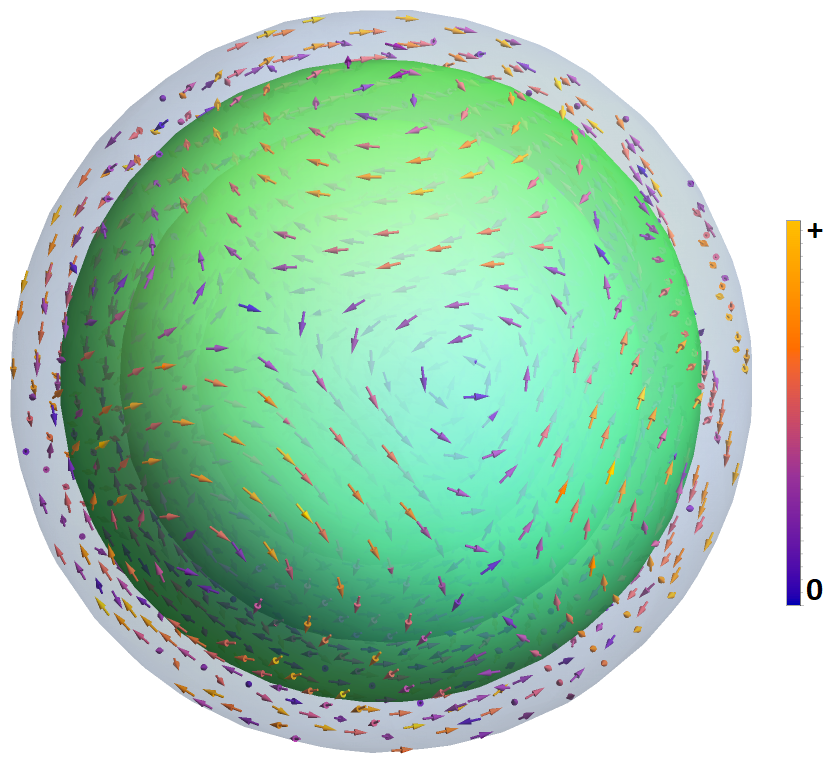}
    \end{subfigure}

    \caption{\emph{Left}: Symplectic foliation of \ec{\Pi_{\mathfrak{so}(3)}} in (\ref{EcPiSO3}). \emph{Right}: Modular vector field of $\Pi$ in (\ref{EcPiNoModSO3}) relative to the Euclidean volume form on \ec{\RR{3}}. The color scale indicates the magnitude of the vectors.}
    \label{fig:SO3-NoUnimod}
\end{figure}

     \subsection{Matrix of a Bivector Field} \label{subsec:matrix}

The value of the matrix (field) \ec{[\Pi^{ij}]} of $\Pi$ at a point \ec{p \in M} is defined to be the matrix
    \begin{equation}\label{EcPiMatrix}
        \big[\Pi^{ij}\big]_{p} := \big[\Pi^{ij}(p)\big].
    \end{equation}
Hence, we just need to know the value of the coefficients of $\Pi$ in (\ref{EcPiCoord}) at $p$.

The function \textsf{num\_bivector\_to\_matrix} evaluates the matrix of a (Poisson) bivector field on a mesh in \ec{\RR{m}}.
\begin{algorithm}[H]
    \captionsetup{justification=centering}
    \caption{\ \textsf{num\_bivector\_to\_matrix}(\emph{bivector, mesh})} \label{AlgNumMatrixBivector}
        \rule{\textwidth}{0.4pt}
    \Input{a (Poisson) bivector field and a mesh}
    \Output{evaluation of the matrix of the (Poisson) bivector field at each point of the mesh \vspace{0.25cm}}
        \rule{\textwidth}{0.4pt}
    \begin{algorithmic}[1] 
        \Procedure{}{}
            \State $m$ $\gets$ dimension of the manifold
            \State \bluecolor{bivector} $\gets$ a variable encoding the (Poisson) bivector field
            \State \bluecolor{mesh} $\gets$ a $(k,m)$ array encoding the mesh
                \CommentNew{$k$: number of points in the mesh}
            \State \bluecolor{variable\_1} $\gets$ a variable encoding the matrix of the (Poisson) bivector field
            \State \textsc{Transform} each \bluecolor{variable\_1} item into a function that allows a numerical evaluation
            \State \rreturn{an array containing the evaluation of \bluecolor{variable\_1} at each $m$--array of \bluecolor{mesh}}
        \EndProcedure
    \end{algorithmic}
\end{algorithm}

\begin{example}
Consider the Lie--Poisson bivector field on \ec{\RR{3}_{x}}
    \begin{equation}\label{EcPiSL2}
        \Pi_{\mathfrak{sl}(2)} =
        -x_3\frac{\partial}{\partial x_1}\wedge \frac{\partial}{\partial x_2} -
        x_2\frac{\partial}{\partial x_1}\wedge \frac{\partial}{\partial x_3} +
        x_1 \frac{\partial}{\partial x_2}\wedge \frac{\partial}{\partial x_3},
    \end{equation}
associated to the 3--dimensional Lie algebra \ec{\mathfrak{sl}(2)}, used in the classification of rigid motions \cite{GaugeClass}, and in other mechanical systems \cite{Naka, Bullo}. To evaluate the matrix of \ec{\Pi_{\mathfrak{sl}(2)}} at points of \ec{Q^{3}} in (\ref{EcCorners}) we compute:
\begin{tcolorbox}[arc=0mm, boxsep=0mm, skin=bicolor, colback=pink!15, colframe=blue!20, colbacklower=blue!0, breakable, halign=left]
\LineCode \textsf{npg3 = NumPoissonGeometry(3)}
    \hspace*{\fill} \CommentCode{\textsf{NumPoissonGeometry} instance} \\
\LineCode \textsf{P\_sl2 = \big\{(1, 2): `-x3', (1, 3): `-x2', (2, 3): `x1'\big\}}
    \hspace*{\fill} \CommentCode{dictionary for \ec{\Pi_{\mathfrak{sl}(2)}} in (\ref{EcPiSL2}) according to (\ref{EcMultivectorDic})} \\
\LineCode \textsf{npg3.num\_bivector\_to\_matrix(P\_sl2, Qmesh, pt\_output=True)} \newline
    \hspace*{\fill} \CommentCode{run \textsf{num\_bivector\_to\_matrix} function with \textsf{pt\_output} flag}
    \tcblower
\resizebox{\textwidth}{!}{
\textsf{tensor([\,\parbox[t]{\linewidth}{[[0.,  0.,  0.], [0.,  0.,  0.], [\phantom{-}0., \phantom{-}0.,  0.]], \
        [[0.,  1.,  0.], [-1.,  0.,  0.], [\phantom{-}0., \phantom{-}0.,  0.]], \newline
        [[0.,  0.,  1.], [0.,  0.,  0.], [-1., \phantom{-}0.,  0.]], \
        [[0.,  1.,  1.], [-1.,  0.,  0.], [-1., \phantom{-}0.,  0.]], \newline
        [[0.,  0.,  0.], [0.,  0.,  1.], [\phantom{-}0., -1.,  0.]], \
        [[0.,  1.,  0.], [-1.,  0.,  1.], [\phantom{-}0., -1.,  0.]], \newline
        [[0.,  0.,  1.], [0.,  0.,  1.], [-1., -1.,  0.]], \
        [[0.,  1.,  1.], [-1.,  0.,  1.], [-1., -1.,  0.]]\,], \newline
        dtype=torch.float64)}}
}
\end{tcolorbox}
\end{example}

    \subsection{Hamiltonian Vector Fields} \label{subsec:hamiltonian}

The Hamiltonian vector field \ec{X_{h} := \ii_{\dd{h}}\Pi} of a scalar function $h$, and relative to a Poisson bivector field $\Pi$ on $M$, at a point in $M$ can be determined by the following (coordinate) formula:
    \begin{equation}\label{EcPiHam}
        X_{h} |_{p} = -\left[\Pi^{ij}(p)\right]\left[\tfrac{\partial h}{\ \partial x^{k}}(p)\right], \quad p \in M
    \end{equation}
Here, \ec{[\Pi^{ij}]} is the matrix of $\Pi$ in (\ref{EcPiMatrix}) and \ec{[{\partial h}/{\partial x^{k}}]} is the gradient vector (field) of $h$, for \ec{k=1,\ldots,m}.

The function \textsf{num\_hamiltonian\_vf} evaluates a Hamiltonian vector field on a mesh in \ec{\RR{m}}.

\begin{algorithm}[H]
    \captionsetup{justification=centering}
    \caption{\ \textsf{num\_hamiltonian\_vf}(\emph{bivector, ham\_function, mesh})} \label{AlgNumHamVF}
        \rule{\textwidth}{0.4pt}
    \Input{a Poisson bivector field $\Pi$, a scalar function $h$ and a mesh}
    \Output{evaluation of the Hamiltonian vector field of $h$ respect to $\Pi$ at each point of the mesh \vspace{0.25cm}}
        \rule{\textwidth}{0.4pt}
    \begin{algorithmic}[1] 
        \Procedure{}{}
            \State $m$ $\gets$ dimension of the manifold
            \State \bluecolor{bivector} $\gets$ a variable encoding the Poisson bivector field $\Pi$
            \State \bluecolor{ham\_function} $\gets$ an expression representing the scalar function $h$
            \State \bluecolor{mesh} $\gets$ a $(k,m)$ array encoding the mesh
                \CommentNew{$k$: number of points in the mesh}
            \State \bluecolor{variable\_1} $\gets$ the (symbolic) gradient vector field of \bluecolor{ham\_function}
            \State \textsc{Transform} each \bluecolor{variable\_1} item into a function that allows a numerical evaluation
            \State \bluecolor{variable\_2} $\gets$ an array containing the evaluation of \bluecolor{variable\_1} at each $m$--array of \bluecolor{mesh}
            \State \bluecolor{variable\_3} $\gets$ \textsf{num\_bivector\_to\_matrix}(\bluecolor{bivector}, \bluecolor{mesh})
                \CommentNew{see Algorithm \ref{AlgNumMatrixBivector}}
            \State \bluecolor{variable\_4} $\gets$ an empty container
            \For {$0 \leq i < k$}
                \State \bluecolor{variable\_4}[$i$] $\gets$ $(-1)$ $*$ \bluecolor{variable\_3}[$i$] $*$ \bluecolor{variable\_2}[$i$]
                    \CommentNew{matrix--vector product}
            \EndFor
            \State \rreturn{\bluecolor{variable\_4}}
        \EndProcedure
    \end{algorithmic}
\end{algorithm}

\begin{example}
Consider the Hamiltonian vector field on \ec{\RR{6}_{x}}, that arises in a particular case of the three body problem \cite{Newton},
     \begin{multline*}
        X_{h} = - x_4\frac{\partial}{\partial{x_1}} - x_5\frac{\partial}{\partial{x_2}} - x_6\frac{\partial}{\partial{x_3}}
        + \left[ \tfrac{1}{(x_1-x_2)|x_1-x_2|} + \tfrac{1}{(x_1-x_3)|x_1-x_3|} \right]\frac{\partial}{\partial{x_4}} \\
        + \left[ \tfrac{1}{(x_1-x_2)|x_1-x_2|} + \tfrac{1}{(x_2-x_3)|x_2-x_3|} \right]\frac{\partial}{\partial{x_5}}
        - \left[ \tfrac{1}{(x_1-x_3)|x_1-x_3|} + \tfrac{1}{(x_2-x_3)|x_2-x_3|} \right]\frac{\partial}{\partial{x_6}},
     \end{multline*}
with Hamiltonian function
    \begin{equation}\label{EcHamR6}
        h = \frac{1}{x_{1}-x_{2}} + \frac{1}{x_{1}-x_{3}} + \frac{1}{x_{2}-x_{3}} + \frac{x_{4}^{2} + x_{5}^{2} + x_{6}^{2}}{2},
    \end{equation}
and relative to the canonical Poisson bivector field on \ec{\RR{6}_{x}}
    \begin{equation}\label{EcPiCanR6}
        \Pi = \frac{\partial}{\partial x_1}\wedge \frac{\partial}{\partial x_4} + \frac{\partial}{\partial x_2}\wedge \frac{\partial}{\partial x_5} + \frac{\partial}{\partial x_3}\wedge \frac{\partial}{\partial x_6}.
    \end{equation}
To evaluate \ec{X_{h}} avoiding singularities we use the \textsf{mesh}, for \ec{Q^{3}} being as in (\ref{EcCorners}),
    \begin{equation*}
        \{-2, -1\} \times \{0, 1\} \times \{2, 3\} \times Q^{3}, \quad \RR{6} = \RR{3}_{\scriptscriptstyle (x_1,x_2,x_3)} \times \RR{3}_{\scriptscriptstyle (x_4,x_5,x_6)}:
    \end{equation*}

\begin{tcolorbox}[arc=0mm, boxsep=0mm, skin=bicolor, colback=pink!15, colframe=blue!20, colbacklower=blue!0, breakable, halign=left]
\LineCode \textsf{npg6 = NumPoissonGeometry(6)}
    \hspace*{\fill} \CommentCode{\textsf{NumPoissonGeometry} instance} \\
\LineCode \textsf{P = \{(1, 4): 1, (2, 5): 1, (3, 6): 1\}}
    \hspace*{\fill} \CommentCode{dictionary for $\Pi$ in (\ref{EcPiCanR6}) according to (\ref{EcMultivectorDic})} \\
\LineCode \textsf{h = \parbox[t]{0.94\linewidth}{`1 / (x1 - x2) + 1 / (x1 - x3) + 1 / (x2 - x3) + (x4**2 + x5**2 + x6**2) / 2'}} \hspace*{\fill} \CommentCode{string expression for $h$ in (\ref{EcHamR6})} \\
\LineCode \textsf{npg6.num\_hamiltonian\_vf(P, h, mesh, pt\_output=True)} \newpage
    \hspace*{\fill} \CommentCode{run \textsf{num\_hamiltonian\_vf} function with \textsf{pt\_output} flag}
    \tcblower
\textsf{tensor([\,\parbox[t]{\linewidth}{[[\phantom{-}0.0000], [\phantom{-}0.0000], [\phantom{-}0.0000], [-0.3125], [0.0000], [0.3125]], \newline
        $\ldots$, \newline
    [[-1.0000], [-1.0000], [\phantom{-}0.0000],  [-0.3125], [0.0000], [0.3125]], \newline
    [[-1.0000], [-1.0000], [-1.0000], [-0.3125], [0.0000],  [0.3125]]\,], \newline
    dtype=torch.float64)}}
\end{tcolorbox}
\end{example}

    \subsection{Poisson Brackets} \label{subsec:bracket}

The Poisson bracket \ec{\{f,g\}_{\Pi}} of two scalar functions $f$ and $g$, induced by a Poisson bivector field $\Pi$ on $M$, at a point in $M$ can be calculated in coordinates following formulas in \cite[Section 1.2]{Dufour}:
    \begin{equation*}
        \{f,g\}_{\Pi}(p) = -\big[\tfrac{\partial g}{\partial x^{k}}(p)\big]^{\top}\left[\Pi^{ij}(p)\right]\big[\tfrac{\partial f}{\partial x^{l}}(p)\big], \quad p \in M
    \end{equation*}
Here, \ec{[\Pi^{ij}]} is the matrix of $\Pi$ in (\ref{EcPiMatrix}), \ec{[{\partial f}/{\partial x^{l}}]} and \ec{[{\partial g}/{\partial x^{k}}]} are the gradient vector (fields) of $f$ and $g$, in that order, for \ec{k,l = 1, \ldots, m}.

The function \textsf{num\_poisson\_bracket} evaluates the Poisson bracket of two scalar functions on a mesh in \ec{\RR{m}}.

\begin{algorithm}[H]
    \captionsetup{justification=centering}
    \caption{\ \textsf{num\_poisson\_bracket}(\emph{bivector, function\_1, function\_2, mesh})} \label{AlgNumPoissonBracket}
        \rule{\textwidth}{0.4pt}
    \Input{a Poisson bivector field $\Pi$, two scalar functions $f,g$ and a mesh}
    \Output{evaluation of the Poisson bracket of $f$ and $g$ induced by $\Pi$ at each point of the mesh \vspace{0.25cm}}
        \rule{\textwidth}{0.4pt}
    \begin{algorithmic}[1] 
        \Procedure{}{}
            \State $m$ $\gets$ dimension of the manifold
            \State \bluecolor{bivector} $\gets$ a variable encoding the Poisson bivector field $\Pi$
            \State \bluecolor{function\_1}, \bluecolor{function\_2} $\gets$ expressions representing the functions $f$ and $g$
            \State \bluecolor{mesh} $\gets$ a $(k,m)$ array encoding the mesh
                \CommentNew{$k$: number of points in the mesh}
            \If {\bluecolor{function\_1} $-$ \bluecolor{function\_2} $== 0$}
                \State \textbf{return} 0
                    \CommentNew{if $f=g$, then its Poisson bracket is zero}
            \Else
                \State \bluecolor{variable\_1} $\gets$ the (symbolic) gradient vector field of \bluecolor{function\_2}
                \State \textsc{Transform} each \bluecolor{variable\_1} item into a function that allows a numerical evaluation
                \State \bluecolor{variable\_2} $\gets$ an array containing the evaluation of \bluecolor{variable\_1} at each $m$--array of \bluecolor{mesh}
                \State \bluecolor{variable\_3} $\gets$ \textsf{num\_hamiltonian\_vf}(\bluecolor{bivector}, \bluecolor{function\_1}, \bluecolor{mesh})
                \Statex    \CommentNew{see Algorithm \ref{AlgNumHamVF}}
                \State \bluecolor{variable\_4} $\gets$ an empty container
                    \breakalg
                \For {$0 \leq i < k$}
                    \State \bluecolor{variable\_4}[$i$] $\gets$ \bluecolor{variable\_2}[$i$] $*$ \bluecolor{variable\_3}[$i$]
                        \CommentNew{scalar vector product}
                \EndFor
                \State \rreturn{\bluecolor{variable\_4}}
            \EndIf
        \EndProcedure
    \end{algorithmic}
\end{algorithm}

\begin{example}
Consider the Poisson bivector field on \ec{\RR{6}_{x}}, obtained as a deformation of an almost Poisson structure analyzed in relation to plasma \cite{Plasma} (see, also \cite{SeveraTwist}),
    \begin{equation}\label{EcTwist}
        \Pi = \frac{\partial}{\partial x_1}\wedge \frac{\partial}{\partial x_4} + \frac{\partial}{\partial x_2}\wedge \frac{\partial}{\partial x_5} + \frac{\partial}{\partial x_3}\wedge \frac{\partial}{\partial x_6} + x_{2}^{2}\frac{\partial}{\partial x_5}\wedge \frac{\partial}{\partial x_6}.
    \end{equation}
Observe that the Poisson bracket \ec{\{x_{6}, x_{5}\}_{\Pi} = -1}, at points \ec{x \in \RR{6}} such that \ec{|x_{2}|=1}. We can check this fact using random \textsf{meshes} of the form
    \begin{equation*}
        \{a_{1}, b_{1}\} \times \{1\} \times \{a_{2}, b_{2}\} \times \cdots \times \{a_{5}, b_{5}\}, \qquad a_{i},b_{i} \in \operatorname{random}{[0,1)}.
    \end{equation*}
Here, random samples are taken from a uniform distribution over the interval \ec{[0,1)}:

\begin{tcolorbox}[arc=0mm, boxsep=0mm, skin=bicolor, colback=pink!15, colframe=blue!20, colbacklower=blue!0, breakable, halign=left]
\LineCode \textsf{npg6 = NumPoissonGeometry(6)}
    \hspace*{\fill} \CommentCode{\textsf{NumPoissonGeometry} instance} \\
\LineCode \textsf{P = \{(1, 4): 1, (2, 5): 1, (3, 6): 1, (5, 6): `x2**2'\}}
    \hspace*{\fill} \CommentCode{dictionary for $\Pi$ in (\ref{EcTwist}) according to (\ref{EcMultivectorDic})} \\
\LineCode \textsf{f, g = `x6', `x5'}
    \hspace*{\fill} \CommentCode{string expressions for canonical coordinates $x_{6}$ and $x_{5}$, in that order} \\
\LineCode \textsf{npg6.num\_poisson\_bracket(P, f, g, meshes, pt\_output=True)} \newline
    \hspace*{\fill} \CommentCode{run the \textsf{num\_poisson\_bracket} function with \textsf{pt\_output} flag}
    \tcblower
\resizebox{\textwidth}{!}{
\textsf{tensor(\parbox[t]{\linewidth}{-1.0, -1.0, -1.0, -1.0, -1.0, -1.0, -1.0, -1.0, -1.0, -1.0, -1.0, -1.0, -1.0, -1.0, -1.0, -1.0, -1.0, -1.0, -1.0, -1.0, -1.0, -1.0, -1.0, -1.0, -1.0, -1.0, -1.0, -1.0, -1.0, -1.0, -1.0, -1.0, -1.0, -1.0, -1.0, -1.0, -1.0, -1.0, -1.0, -1.0, -1.0, -1.0, -1.0, -1.0, -1.0, -1.0, -1.0, -1.0, -1.0, -1.0, -1.0, -1.0, -1.0, -1.0, -1.0, -1.0, -1.0, -1.0, -1.0, -1.0, -1.0, -1.0, -1.0, -1.0, dtype=torch.float64)}}
    }
\end{tcolorbox}
\end{example}

Observe that we use a probabilistic method to verify a particular property of a Poisson bracket. These methods can be used to determine other characteristics, for example, singular points.

    \subsection{Sharp Morphism} \label{subsec:sharp}

In the context of Lie algebroids, the vector bundle morphism \ec{\Pi^{\sharp}:\alpha \mapsto \ii_{\alpha}\Pi}, induced by a Poisson bivector field $\Pi$ on $M$, is the anchor map of the Poisson Lie algebroid corresponding to $\Pi$ \cite{FuchAlgbd}.

Similarly to (\ref{EcPiHam}), we can evaluate the image \ec{\Pi^{\sharp}(\alpha)} of a differential 1--form $\alpha$ at a point in $M$ as follows:
    \begin{equation}\label{EcPiSharp}
        \Pi^{\sharp}(\alpha)\big|_{p} = -\left[\Pi^{ij}(p)\right]\big[\alpha_{k}(p)\big], \quad p \in M
    \end{equation}
Here, \ec{[\Pi^{ij}]} is the matrix of $\Pi$ in (\ref{EcPiMatrix}) and \ec{[\alpha_{k}]} is the coefficient vector of \ec{\alpha = \alpha_{k}\dd{x^{k}}}, for \ec{k = 1, \ldots, m}.

The function \textsf{num\_sharp\_morphism} evaluates a vector field \ec{\Pi^{\sharp}(\alpha)} on a mesh in \ec{\RR{m}}.

\begin{algorithm}[H]
    \captionsetup{justification=centering}
    \caption{\ \textsf{num\_sharp\_morphism}(\emph{bivector, one\_form, mesh})} \label{AlgNumSharp}
        \rule{\textwidth}{0.4pt}
    \Input{a Poisson bivector field $\Pi$, a differential 1--form $\alpha$ and a mesh}
    \Output{evaluation of the vector field $\Pi^{\sharp}(\alpha)$ at each point of the mesh}
        \rule{\textwidth}{0.4pt}
    \begin{algorithmic}[1] 
        \Procedure{}{}
            \State $m$ $\gets$ dimension of the manifold
            \State \bluecolor{bivector} $\gets$ a variable encoding the Poisson bivector field $\Pi$
            \State \bluecolor{one\_form} $\gets$ a variable encoding the differential 1--form \ec{\alpha = \alpha_{1}\dd{x^{1}} + \cdots + \alpha_{m}\dd{x^{m}}}
            \State \bluecolor{mesh} $\gets$ a $(k,m)$ array encoding the mesh
                \CommentNew{$k$: number of points in the mesh}
            \State \bluecolor{variable\_1} $\gets$ a container with items $(\alpha_{1},\ldots,\alpha_{m})$
            \State \textsc{Transform} each \bluecolor{variable\_1} item into a function that allows a numerical evaluation
            \State \bluecolor{variable\_2} $\gets$ an array containing the evaluation of \bluecolor{variable\_1} at each $m$--array of \bluecolor{mesh}
            \State \bluecolor{variable\_3} $\gets$ \textsf{num\_bivector\_to\_matrix}(\bluecolor{bivector}, \bluecolor{mesh})
                \CommentNew{see Algorithm \ref{AlgNumMatrixBivector}}
            \State \bluecolor{variable\_4} $\gets$ an empty container
            \For {$0 \leq i < k$}
                \State \bluecolor{variable\_4}[$i$] $\gets$ $(-1)$ $*$ \bluecolor{variable\_3}[$i$] $*$ \bluecolor{variable\_2}[$i$]
                    \CommentNew{matrix--vector product}
            \EndFor
                \State \rreturn{\bluecolor{variable\_4}}
        \EndProcedure
    \end{algorithmic}
\end{algorithm}

\begin{example}
Consider the Casimir function \ec{K = ({1}/{2})(x_{1}^{2} + x_{2}^{2} + x_{3}^{2})} of the Poisson bivector field \ec{\Pi_{\mathfrak{so}(3)}} in (\ref{EcPiSO3}) \cite{Damianou,Newton}. By definition, the exterior derivative of $K$,
    \begin{equation}\label{EcDKsharp}
        \dd{K} = x_{1}\dd{x_{1}} + x_{2}\dd{x_{2}} + x_{3}\dd{x_{3}},
    \end{equation}
belongs to the kernel of \ec{\Pi_{\mathfrak{so}(3)}}, that is, \ec{\Pi_{\mathfrak{so}(3)}^{\sharp}(\dd{K}) = 0}. We can check this fact by using random meshes:

\begin{tcolorbox}[arc=0mm, boxsep=0mm, skin=bicolor, colback=pink!15, colframe=blue!20, colbacklower=blue!0, breakable, halign=left]
\LineCode \textsf{npg3 = NumPoissonGeometry(3)}
    \hspace*{\fill} \CommentCode{\textsf{NumPoissonGeometry} instance} \\
\LineCode \textsf{P\_so3 = \big\{(1, 2): `x3', (1, 3): `-x2', (2, 3): `x1'\big\}}
    \hspace*{\fill} \CommentCode{dictionary for $\Pi_{\mathfrak{so}(3)}$ in (\ref{EcPiSO3}) according to (\ref{EcMultivectorDic})} \\
\LineCode \textsf{dK = \{(1,): `x1', (2,): `x2', (3,): `x3'\}}
    \hspace*{\fill} \CommentCode{dictionary for $\dd{K}$ in (\ref{EcDKsharp}) according to (\ref{EcMultivectorDic})} \\
\LineCode \textsf{mesh = numpy.random.rand(10**6, 3)} \\
    \hspace*{\fill} \CommentCode{$(10^{6}, 3)$ \textsf{NumPy} array with random samples from a uniform distribution over [0,1)} \\
\LineCode \textsf{npg6.num\_sharp\_morphism(P\_so3, dK, mesh, pt\_output=True)} \newline
    \hspace*{\fill} \CommentCode{run \textsf{num\_sharp\_morphism} function with \textsf{pt\_output} flag}
    \tcblower
\textsf{tensor([\,\parbox[t]{\linewidth}{[[0.], [0.], [0.]], \ [[0.], [0.], [0.]], \ [[0.], [0.], [0.]], \ $\ldots$, \newline
    [[0.], [0.], [0.]], \ [[0.], [0.], [0.]], \ [[0.], [0.], [0.]]\,], dtype=torch.float64)}}
\end{tcolorbox}
\end{example}

    \subsection{Coboundary Operator} \label{subsec:coboundary}

The coboundary operator $\delta_{\Pi}:\Gamma\wedge^{\bullet}\T{M} \rightarrow \Gamma\wedge^{\bullet + 1}\T{M}$ induced by a Poisson bivector field $\Pi$ on $M$ \cite{Lich-77} is defined by
     \begin{equation*}
          \delta_{\Pi}(A) := \cSch{\Pi,A}, \quad A \in \Gamma\wedge\T{M}.
     \end{equation*}
Here, \ec{\Gamma\wedge^{\bullet}\T{M}} denotes the module of multivector fields on $M$.

The function \textsf{num\_coboundary\_operator} evaluates the image under \ec{\delta_{\Pi}} of an arbitrary multivector field on a mesh in \ec{\RR{m}}.

\begin{algorithm}[H]
    \captionsetup{justification=centering}
    \caption{\ \textsf{num\_coboundary\_operator}(\emph{bivector, multivector, mesh})} \label{AlgNumCoboundary}
        \rule{\textwidth}{0.4pt}
    \Input{a Poisson bivector field $\Pi$, a multivector field $A$ and a mesh}
    \Output{evaluation of the multivector field $\cSch{\Pi, A}$ at each point of the mesh}
        \rule{\textwidth}{0.4pt}
    \begin{algorithmic}[1] 
        \Procedure{}{}
            \State $m$ $\gets$ dimension of the manifold
            \State \bluecolor{bivector} $\gets$ a variable encoding the Poisson bivector field $\Pi$
            \State \bluecolor{multivector} $\gets$ a variable encoding the multivector field $A$
            \State \bluecolor{mesh} $\gets$ a $(k,m)$ array encoding the mesh
                \CommentNew{$k$: number of points in the mesh}
            \State \bluecolor{variable\_1} $\gets$ a variable encoding the image of \bluecolor{multivector} under the coboundary operator induced by \bluecolor{bivector}
            \Statex \CommentNew{We use the \textsf{lichnerowicz\_poisson\_operator} method of \textsf{PoissonGeometry}}
            \State \textsc{Transform} each \bluecolor{variable\_1} item into a function that allows a numerical evaluation
            \State \rreturn{an array containing the evaluation of \bluecolor{variable\_1} at each $m$--array of \bluecolor{mesh}}
        \EndProcedure
    \end{algorithmic}
\end{algorithm}

\begin{example}
The vector field on \ec{\RR{3}_{x} \setminus \{x_{3}\text{--axis}\}}
    \begin{multline}\label{EcWCobound}
        W = \frac{x_1x_{3}e^{{-1}/{(x_1^2 + x_2^2 - x_3^2)^2}}}{x_1^2 + x_2^2}\frac{\partial}{\partial{x_{1}}} +
        \frac{x_2x_{3}e^{{-1}/{(x_1^2 + x_2^2 - x_3^2)^2}}}{x_1^2 + x_2^2}\frac{\partial}{\partial{x_{2}}} \\
        + e^{{-1}/{(x_1^2 + x_2^2 - x_3^2)^2}}\frac{\partial}{\partial{x_{3}}}
    \end{multline}
arises in the study of the first cohomology group of the Poisson bivector field \ec{\Pi_{\mathfrak{sl}(2)}} in (\ref{EcPiSL2}) to construct a 1--cocycle that is not Hamiltonian \cite{Naka, MarcutSl2}. To check the cocycle property of $W$ under the coboundary operator induced by \ec{\Pi_{\mathfrak{sl}(2)}}, we evaluate the (image) bivector field \ec{\cSch{\Pi_{\mathfrak{sl}(2)},W}} on random meshes:

\begin{tcolorbox}[arc=0mm, boxsep=0mm, skin=bicolor, colback=pink!15, colframe=blue!20, colbacklower=blue!0, breakable, halign=left]
\LineCode \textsf{npg3 = NumPoissonGeometry(3)}
    \hspace*{\fill} \CommentCode{\textsf{NumPoissonGeometry} instance} \\
\LineCode \textsf{P\_sl2 = \big\{(1, 2): `-x3', (1, 3): `-x2', (2, 3): `x1'\big\}}
    \hspace*{\fill} \CommentCode{dictionary for $\Pi_{\mathfrak{sl}(2)}$ in (\ref{EcPiSL2}) according to (\ref{EcMultivectorDic})} \\
\LineCode \textsf{W = \{\parbox[t]{0.95\linewidth}{(1,): `x1 * x3 * exp(-1 / (x1**2 + x2**2 - x3**2)**2) / (x1**2 + x2**2)', \newline
(2,): `x2 * x3 * exp(-1 / (x1**2 + x2**2 - x3**2)**2) / (x1**2 + x2**2)', \newline
(3,): `exp(-1 / (x1**2 + x2**2 - x3**2)**2)'\}}}
    \hspace*{\fill} \CommentCode{dictionary for $W$ in (\ref{EcWCobound}) according to (\ref{EcMultivectorDic})} \\
\LineCode \textsf{mesh = numpy.random.rand(10**6, 3)} \\
    \hspace*{\fill} \CommentCode{$(10^{6}, 3)$ \textsf{NumPy} array with random samples from a uniform distribution over [0,1)} \\
\LineCode \textsf{npg3.num\_coboundary\_operator(P\_sl2, W, mesh, pt\_output=True)} \newline
    \hspace*{\fill} \CommentCode{run \textsf{num\_coboundary\_operator} function with \textsf{pt\_output} flag}
    \tcblower
\textsf{tensor([\,\parbox[t]{\linewidth}{[[0.,  0.,  0.], [0.,  0.,  0.], [0., 0.,  0.]], \
        [[0.,  0.,  0.], [0.,  0.,  0.], [0.,  0.,  0.]], \newline
                        $\ldots$, \newline
        [[0.,  0.,  0.], [0.,  0.,  0.], [0.,  0.,  0.]], \
        [[0.,  0.,  0.], [0.,  0.,  0.], [0.,  0.,  0.]], \newline
        [[0.,  0.,  0.], [0.,  0.,  0.], [0.,  0.,  0.]], \
        [[0.,  0.,  0.], [0.,  0.,  0.], [0.,  0.,  0.]]\,], \newline
        dtype=torch.float64)}}
\end{tcolorbox}
\end{example}

\begin{figure}[H]
\centering
    \begin{subfigure}{0.4\textwidth}
        \includegraphics[width=\textwidth]{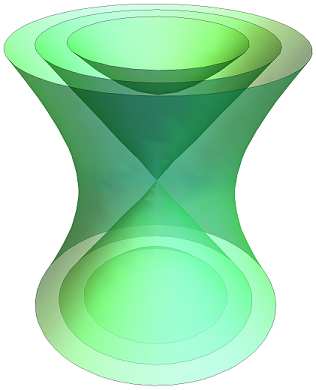}
    \end{subfigure}
        \hspace{0.1cm}
    \begin{subfigure}{0.4\textwidth}
        \includegraphics[width=\textwidth]{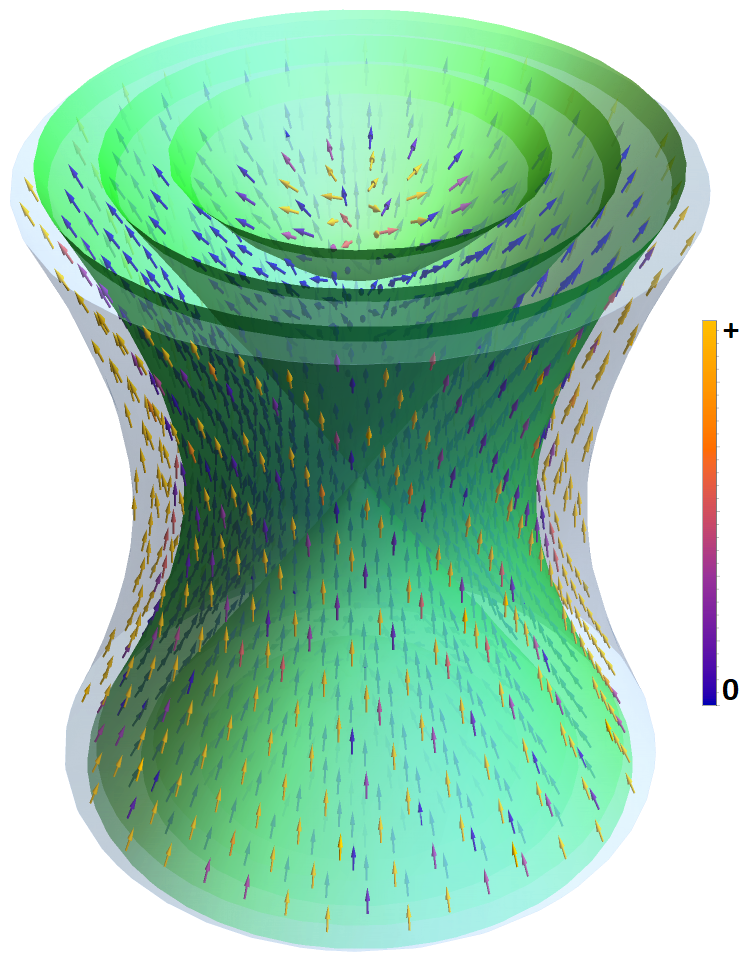}
    \end{subfigure}

    \caption{\emph{Left}: Symplectic foliation of \ec{\Pi_{\mathfrak{sl}(2)}} in (\ref{EcPiSL2}). \emph{Right}: Vector field W in (\ref{EcWCobound}), tangent to the symplectic foliation of \ec{\Pi_{\mathfrak{sl}(2)}} on \ec{\RR{3}_{x} \setminus \{x_{3}\text{--axis}\}}. The color scale indicates the magnitude of the vectors.}
    \label{fig:SL2-TangentVF}
\end{figure}

The characteristic foliation of \ec{\Pi_{\mathfrak{sl}(2)}} can be described using the Casimir function \ec{2K=x_{1}^{2} + x_{2}^{2} - x_{3}^{2}} (see, Figure \ref{fig:SL2-TangentVF}). Observe that $W$ is orthogonal to the gradient vector field of $K$, \ec{x_{1}\partial/\partial{x_{1}} + x_{2}\partial/\partial{x_{2}} - x_{3}\partial/\partial{x_{3}}}. Hence, $W$ is tangent to the symplectic foliation of \ec{\Pi_{\mathfrak{sl}(2)}} (see, Figure \ref{fig:SL2-TangentVF}). However, $W$ can be extended to a smooth vector field on all of \ec{\RR{3}} that is not Hamiltonian with respect to \ec{\Pi_{\mathfrak{sl}(2)}} \cite{Naka}. Therefore the first Poisson cohomology group of \ec{\Pi_{\mathfrak{sl}(2)}} is non--trivial.

    \subsection{Modular Vector Field} \label{subsec:modular}

The modular vector field of an orientable Poisson manifold \ec{(M,\Pi)} is an infinitesimal automorphism of $\Pi$ determined by the choice of a volume form $\Omega$ \cite{WeModular}, and it is defined by the linear map
    \begin{equation*}
        Z_{\Pi, \Omega}: h \longmapsto \mathrm{div}_{\Omega}X_{h}, \quad h \in \Cinf{M}.
    \end{equation*}
Here, \ec{\mathrm{div}_{\Omega}X_{h} \in \Cinf{M}} denotes the divergence of the Hamiltonian vector field \ec{X_{h}} (\ref{EcPiHam}) with respect to $\Omega$.

The modular vector field measures how far Hamiltonian flows are from preserving a given volume form \cite{WeModular}. In the regular case, the existence of a volume form that remains invariant under every Hamiltonian vector field only depends on the characteristic (symplectic) foliation of the Poisson manifold, rather than the leaf--wise symplectic form \cite{Reeb2}. In particular, for regular codimension--one symplectic foliations, a characteristic class controls the existence of such a volume form \cite{Miranda}.

The function \textsf{num\_modular\_vf} evaluates the modular vector field of a Poisson bivector field on a mesh in \ec{\RR{m}}.

\begin{algorithm}[H]
    \captionsetup{justification=centering}
    \caption{\ \textsf{num\_modular\_vf}(\emph{bivector, function, mesh})} \label{AlgNumModularVF}
        \rule{\textwidth}{0.4pt}
    \Input{a Poisson bivector field $\Pi$, a non--zero scalar function $f_{0}$ and a mesh}
    \Output{evaluation of the modular vector field of $\Pi$ relative to the volume form $f_{0}\Omega_{0}$ at each point of the mesh. Where $\Omega_{0}$ is the Euclidean volume form on $\RR{m}$ \vspace{0.25cm}}
        \rule{\textwidth}{0.4pt}
    \begin{algorithmic}[1] 
        \Procedure{}{}
            \State $m$ $\gets$ dimension of the manifold
            \State \bluecolor{bivector} $\gets$ a variable encoding the Poisson bivector field $\Pi$
            \State \bluecolor{function} $\gets$ an expression representing the scalar function $f_{0}$
            \State \bluecolor{mesh} $\gets$ a $(k,m)$ array encoding the mesh
                \CommentNew{$k$: number of points in the mesh}
             \State \bluecolor{variable\_1} $\gets$ a variable encoding the modular vector field of \bluecolor{bivector} with respect to the volume form (\bluecolor{function})$\cdot \Omega_{0}$
                \CommentNew{We use the \textsf{modular\_vf} method of \textsf{PoissonGeometry}}
            \State \textsc{Transform} each \bluecolor{variable\_1} item into a function that allows a numerical evaluation
            \State \rreturn{an array containing the evaluation of \bluecolor{variable\_1} at each $m$--array of \bluecolor{mesh}}
        \EndProcedure
    \end{algorithmic}
\end{algorithm}

\begin{example}\label{example:NoModSO3}
The characteristic foliation of the following homogeneous Poisson bivector field on \ec{\RR{3}_{x}} coincides with that of \ec{\Pi_{\mathfrak{so}(3)}} in (\ref{EcPiSO3}):
    \begin{multline}\label{EcPiNoModSO3}
        \Pi =
        \tfrac{1}{4}x_3\big(x_{1}^{4} + x_{2}^{4} + x_{3}^{4}\big)\frac{\partial}{\partial x_1}\wedge \frac{\partial}{\partial x_2} -
        \tfrac{1}{4}x_2\big(x_{1}^{4} + x_{2}^{4} + x_{3}^{4}\big)\frac{\partial}{\partial x_1}\wedge \frac{\partial}{\partial x_3} \\
        + \tfrac{1}{4}x_1\big(x_{1}^{4} + x_{2}^{4} + x_{3}^{4}\big)\frac{\partial}{\partial x_2}\wedge \frac{\partial}{\partial x_3}
    \end{multline}
However, \ec{(\RR{3},\Pi_{\mathfrak{so}(3)})} admits a Hamiltonian--invariant volume form, while $\Pi$ does not \cite{Pedroza}. Consequently, the following vector field,
    \begin{equation*}
        Z_{\Pi} = \big(x_{2}x_{3}^{3} - x_{2}^{3}x_{3}\big)\frac{\partial}{\partial{x_{1}}} +
        \big(x_{3}x_{1}^{3} - x_{3}^{3}x_{1}\big)\frac{\partial}{\partial{x_{2}}} +
        \big(x_{1}x_{2}^{3} - x_{1}^{3}x_{2}\big)\frac{\partial}{\partial{x_{3}}},
    \end{equation*}
which is the modular vector field of $\Pi$ with respect to the Euclidean volume form on \ec{\RR{3}_{x}}, cannot be a Hamiltonian vector field.

We can now use random meshes to check numerically that there are points for which the modular vector field of $\Pi$ with respect to the Euclidean volume form is not zero, which implies that $\Pi$ is not unimodular, as follows:

\begin{tcolorbox}[arc=0mm, boxsep=0mm, skin=bicolor, colback=pink!15, colframe=blue!20, colbacklower=blue!0, breakable, halign=left]
\LineCode \textsf{npg3 = NumPoissonGeometry(3)}
    \hspace*{\fill} \CommentCode{\textsf{NumPoissonGeometry} instance} \newline
\LineCode \textsf{P = \{\parbox[t]{0.85\linewidth}{(1, 2): `1/4*x3*(x1**4 + x2**4 + x3**4)', (1, 3): `- 1/4*x2* (x1**4 + x2**4 + x3**4)', (2, 3): `1/4*x1*(x1**4 + x2**4 + x3**4)'\}}} \hspace*{\fill} \CommentCode{dictionary for $\Pi$ in (\ref{EcPiNoModSO3}) according to (\ref{EcMultivectorDic})} \\
\LineCode \textsf{mesh = numpy.random.rand(10**6, 3)} \\
    \hspace*{\fill} \CommentCode{$(10^{6}, 3)$ \textsf{NumPy} array with random samples from a uniform distribution over [0,1)} \\
\LineCode \textsf{npg3.num\_modular\_vf(P, 1, mesh, pt\_output=True)}
    \hspace*{\fill} \CommentCode{run \textsf{num\_modular\_vf} function with \textsf{pt\_output} flag}
    \tcblower
\textsf{tensor([\,\parbox[t]{\linewidth}{[[ 0.0538], [ 0.0545], [-0.1005]], \ [[ 0.0031], [\phantom{-}0.3838], [-0.0149]], \newline
    [[-0.2559], [-0.0204], [ 0.0575]], \ [[-0.0910], [\phantom{-}0.0093], [ 0.1988]], \newline
        $\ldots$, \newline
    [[-0.0013], [-0.0095], [ 0.0107]], \ [[-0.0467], [-0.0012], [0.0296]]\,], \newline
    dtype=torch.float64)}}
\end{tcolorbox}
\end{example}

Observe that \ec{Z_{\Pi}} is orthogonal to the radial vector field \ec{x_{i}\partial/\partial{x_{i}}} on \ec{\RR{3}_{x}}, for \ec{i=1,2,3}. As the characteristic foliation of $\Pi$ consists of the origin and concentric spheres (see, Figure \ref{fig:SO3-NoUnimod}), \ec{Z_{\Pi}} is tangent to the symplectic leafs of $\Pi$.

    \subsection{Curl Operator} \label{subsec:curl}

On an oriented manifold $M$ with volume form $\Omega$, the \emph{divergence} of an $a$--multivector field $A$ on $M$ \cite{Kozul}, and relative to $\Omega$, is the unique \ec{(a-1)}--multivector field \ec{\mathscr{D}_{\Omega}(A)} on \ec{M} such that
    \begin{equation}\label{EcTrazaDef}
        \ii_{\mathscr{D}_{\Omega}(A)}\Omega = \dd{\ii_{A}}\Omega.
    \end{equation}
The function \textsf{num\_curl\_operator} evaluates the divergence of a multivector field on a mesh in \ec{\RR{m}}. Let \ec{\Omega_{0}} denote the standard volume form on \ec{\RR{m}}:

\begin{algorithm}[H]
    \captionsetup{justification=centering}
    \caption{\ \textsf{num\_curl\_operator}(\emph{multivector, function, mesh})} \label{AlgNumCurl}
        \rule{\textwidth}{0.4pt}
    \Input{a multivector field $A$, a non--zero scalar function $f_{0}$ and a mesh}
    \Output{evaluation of the the divergence of $A$ with respect to the volume form $f_{0}\Omega_{0}$ at each point of the mesh. \vspace{0.25cm}}
        \rule{\textwidth}{0.4pt}
    \begin{algorithmic}[1] 
        \Procedure{}{}
            \State $m$ $\gets$ dimension of the manifold
            \State \bluecolor{multivector} $\gets$ a variable encoding the multivector field $A$
            \State \bluecolor{function} $\gets$ an expression representing the scalar function $f_{0}$
            \State \bluecolor{mesh} $\gets$ a $(k,m)$ array encoding the mesh
                \CommentNew{$k$: number of points in the mesh}
            \State \bluecolor{variable\_1} $\gets$ a variable encoding the divergence of \bluecolor{bivector} with respect to the volume form (\bluecolor{function})$\cdot \Omega_{0}$
                \CommentNew{We use the \textsf{curl\_operator} method of \textsf{PoissonGeometry}}
            \State \textsc{Transform} each \bluecolor{variable\_1} item into a function that allows a numerical evaluation
            \State \rreturn{an array containing the evaluation of \bluecolor{variable\_1} at each $m$--array of \bluecolor{mesh}}
        \EndProcedure
    \end{algorithmic}
\end{algorithm}

\begin{example}
The next Poisson bivector field on \ec{\RR{4}_{x}}, has been applied to the analyze the orbital stability of the Pais--Uhlenbeck oscillator \cite{MVallYu},
    \begin{multline}\label{EcPiPaisU}
        \Pi = 2x_{4}\frac{\partial}{\partial x_1}\wedge \frac{\partial}{\partial x_3}
        + 2x_{3} \frac{\partial}{\partial x_1}\wedge \frac{\partial}{\partial x_4}
        - 2x_{4} \frac{\partial}{\partial x_2}\wedge \frac{\partial}{\partial x_3}
        + 2x_{3} \frac{\partial}{\partial x_2}\wedge \frac{\partial}{\partial x_4} \\
        + (x_{1}-x_{2}) \frac{\partial}{\partial x_3}\wedge \frac{\partial}{\partial x_4}.
    \end{multline}
It is unimodular on (the whole of) \ec{\RR{4}_{x}}, and has trivial modular vector field with respect to the Euclidean volume form. For a fixed volume form, the divergence of a Poisson bivector field coincides with its (minus) modular vector field. Hence, the unimodularity of $\Pi$ in (\ref{EcPiPaisU}) may be numerically verified, at least in a sample, using random meshes. In this example the output is a Pytorch tensor:

\begin{tcolorbox}[arc=0mm, boxsep=0mm, skin=bicolor, colback=pink!15, colframe=blue!20, colbacklower=blue!0, breakable, halign=left]
\LineCode \textsf{npg4 = NumPoissonGeometry(4)}
    \hspace*{\fill} \CommentCode{\textsf{NumPoissonGeometry} instance} \\
\LineCode \textsf{P = \parbox[t]{0.9\linewidth}{\big\{(1, 3): `2*x4', (1, 4): `2*x3', (2, 3): `-2*x4', (2, 4): `2*x3', (3, 4): `x1 - x2'\}}} \hspace*{\fill} \CommentCode{dictionary for $\Pi$ in (\ref{EcPiPaisU}) according to (\ref{EcMultivectorDic})} \\
\LineCode \textsf{mesh = numpy.random.rand(10**6, 4)} \\
    \hspace*{\fill} \CommentCode{$(10^{6}, 4)$ \textsf{NumPy} array with random samples from a uniform distribution over [0, 1)} \\
\LineCode \textsf{npg4.num\_curl\_operator(P, 1, mesh, pt\_output=True)} \newline
    \hspace*{\fill} \CommentCode{run \textsf{num\_curl\_operator} function with \textsf{pt\_output} flag}
    \tcblower
\textsf{tensor([\,\parbox[t]{\linewidth}{[[0.], [0.], [0.], [0.]], \ [[0.], [0.], [0.], [0.]], \ [[0.], [0.], [0.], [0.]], \ $\ldots$, \newline
    [[0.], [0.], [0.], [0.]], \ [[0.], [0.], [0.], [0.]], \ [[0.], [0.], [0.], [0.]]\,], \newline
    dtype=torch.float64)}}
\end{tcolorbox}
\end{example}

    \subsection{Differential 1--Forms Bracket} \label{subsec:onefbracket}

In the context of Lie algebroids, the Koszul bracket of 1--forms is the Lie bracket on the space of sections of a Poisson Lie algebroid \cite{Kozul, Grabowski}.

By definition and (\ref{EcPiSharp}), the Lie bracket of two differential 1--forms $\alpha$ and $\beta$, induced by a Poisson bivector field $\Pi$ on $M$, at a point in $M$ can be determined by the following (coordinate) formula:
    \begin{multline*}
      \{\alpha,\beta\}_{\Pi}(p) = \big[(\dd{\beta})_{kl}(p)\big]\left[\Pi^{ij}(p)\right]\big[\alpha_{r}(p)\big] - \big[(\dd{\alpha})_{kl}(p)\big]\left[\Pi^{ij}(p)\right]\big[\beta_{r}(p)\big] \\
      - \big[\tfrac{\partial}{\partial x^{r}}\big(\left[\beta_{k}\right]^{\top}\left[\Pi^{ij}\right]\left[\alpha_{l}\right]\big)(p)\big], \quad p \in M
    \end{multline*}
Here, \ec{[\Pi^{ij}]} is the matrix of $\Pi$ in (\ref{EcPiMatrix}), \ec{[(\dd{\alpha})_{kl}]} and \ec{[(\dd{\beta})_{kl}]} are the matrix of the differential 2--forms \ec{\dd{\alpha} = 1/2\,(\dd{\alpha})_{kl}\dd{x^{k}} \wedge \dd{x^{l}}} and \ec{\dd{\beta} = 1/2\,(\dd{\beta})_{kl}\dd{x^{k}} \wedge \dd{x^{l}}}, \ec{[\alpha_{r}]} and \ec{[\beta_{r}]} are the coefficient vectors of \ec{\alpha = \alpha_{k}\dd{x^{k}}} and \ec{\beta = \beta_{k}\dd{x^{k}}}, and \ec{[{\partial}/{\partial x^{r}}(\cdot)]} is the gradient vector (field) operator, for \ec{k,l,r = 1, \ldots, m}.

The function \textsf{num\_one\_forms\_bracket} evaluates the differential 1--form \ec{\{\alpha,\beta\}_{\Pi}} on a mesh in \ec{\RR{m}}.

\begin{algorithm}[H]
    \captionsetup{justification=centering}
    \caption{\ \textsf{num\_one\_forms\_bracket}(\emph{bivector, one\_form\_1, one\_form\_2, mesh})} \label{AlgNumOneFormsB}
        \rule{\textwidth}{0.4pt}
    \Input{a Poisson bivector field $\Pi$, two differential 1--forms $\alpha, \beta$ and a mesh}
    \Output{evaluation of the Lie bracket of $\alpha$ and $\beta$ induced by $\Pi$ at each point of the mesh \vspace{0.25cm}}
        \rule{\textwidth}{0.4pt}
    \begin{algorithmic}[1] 
        \Procedure{}{}
            \State $m$ $\gets$ dimension of the manifold
            \State \bluecolor{bivector} $\gets$ a variable encoding the Poisson bivector field $\Pi$
            \State \bluecolor{one\_form\_1} $\gets$ a variable encoding the differential 1--form \ec{\alpha = \alpha_{1}\dd{x^{1}} + \cdots + \alpha_{m}\dd{x^{m}}}
            \State \bluecolor{one\_form\_2} $\gets$ a variable encoding the differential 1--form \ec{\beta = \beta_{1}\dd{x^{1}} + \cdots + \beta_{m}\dd{x^{m}}}
            \State \bluecolor{mesh} $\gets$ a $(k,m)$ array encoding the mesh
                \CommentNew{$k$: number of points in the mesh}
            \State \bluecolor{variable\_1} $\gets$ a container with items $(\alpha_{1},\ldots,\alpha_{m})$
            \State \bluecolor{variable\_2} $\gets$ a container with items $(\beta_{1},\ldots,\beta_{m})$
            \For {$i \in \{1,2\}$}
                \State \bluecolor{variable\_3\_}$i$ $\gets$ the (symbolic) Jacobian matrix of \bluecolor{variable\_}$i$
                \Statex    \CommentNew{\bluecolor{variable\_}$i$ thought of as a vector field}
                \State \textsc{Transform} each \bluecolor{variable\_3\_i} item into a function that allows a numerical evaluation
                \State \bluecolor{variable\_4\_}$i$ $\gets$ an array containing the evaluation of \bluecolor{variable\_3\_}$i$ at each $m$--array of \bluecolor{mesh}
                \State \bluecolor{variable\_5\_}$j$ $\gets$ \textsf{num\_sharp\_morphism}(\bluecolor{bivector}, \bluecolor{one\_form\_}$j$, \bluecolor{mesh})
                \Statex \CommentNew{\ec{i \neq j = 1,2}. See Algorithm \ref{AlgNumMatrixBivector}}
            \EndFor
            \State \bluecolor{variable\_6}, \bluecolor{variable\_7} $\gets$ empty containers
            \For {$0 \leq i < k$}
                \State \bluecolor{variable\_6}[$i$] $\gets$ \big(\bluecolor{variable\_4\_1} $-$ \textsf{transpose}(\bluecolor{variable\_4\_1})\big)[$i$] $*$ \bluecolor{variable\_5\_2}[$i$]
                \State \bluecolor{variable\_7}[$i$] $\gets$ \big(\bluecolor{variable\_4\_2} $-$ \textsf{transpose}(\bluecolor{variable\_4\_2})\big)[$i$] $*$ \bluecolor{variable\_5\_1}[$i$]
            \Statex \CommentNew{matrix--vector products}
            \EndFor
                \breakalg
            \State \bluecolor{variable\_8} $\gets$ a variable encoding the image of \bluecolor{one\_form\_}$j$ under the sharp morphism induced by \bluecolor{bivector}
                \CommentNew{We use the \textsf{sharp\_morphism} method of \textsf{PoissonGeometry}}
            \State \bluecolor{variable\_9} $\gets$ \bluecolor{variable\_8} $*$ \bluecolor{variable\_2}
                \CommentNew{scalar vector product}
            \State \bluecolor{variable\_10} $\gets$ the (symbolic) gradient vector field of \bluecolor{variable\_9}
            \State \textsc{Transform} each \bluecolor{variable\_10} item into a function that allows a numerical evaluation
            \State \bluecolor{variable\_11} $\gets$ an array containing the evaluation of \bluecolor{variable\_10} at each $m$--array of \bluecolor{mesh}
            \State \bluecolor{variable\_12} $\gets$ an empty container
            \For {$0 \leq i < k$}
                \State \bluecolor{variable\_12}[$i$] $\gets$ \bluecolor{variable\_7}[$i$] $-$ \bluecolor{variable\_6}[$i$] $+$ \bluecolor{variable\_11}[$i$]
                    \CommentNew{vector sum}
            \EndFor
            \State \rreturn{\bluecolor{variable\_12}}
        \EndProcedure
    \end{algorithmic}
\end{algorithm}

\begin{example}
By definition, the Lie bracket induced by the Poisson bivector field $\Pi$ in (\ref{EcTwist}) \cite{Plasma}, see also \cite{SeveraTwist}, of the (basic) differential 1--forms \ec{\dd{x_{5}}} and \ec{\dd{x_{6}}} is given by
    \begin{equation*}
        \{\dd{x_{5}}, \dd{x_{6}}\}_{\Pi} = 2x_{2}\dd{x_{2}}.
    \end{equation*}
Hence, \ec{\{\dd{x_{5}}, \dd{x_{6}}\}_{\Pi} = 2\dd{x_{2}}} at points \ec{x \in \RR{6}} such that \ec{x_{2}=1}. This coincides with the following computation involving random \textsf{mesh}es,
    \begin{equation*}
        \{a_{1}, b_{1}\} \times \{1\} \times \{a_{2}, b_{2}\} \times \cdots \times \{a_{5}, b_{5}\}, \quad a_{i},b_{i} \in \operatorname{random}{[0,1]}:
    \end{equation*}

\begin{tcolorbox}[arc=0mm, boxsep=0mm, skin=bicolor, colback=pink!15, colframe=blue!20, colbacklower=blue!0, breakable, halign=left]
\LineCode \textsf{npg6 = NumPoissonGeometry(6)}
    \hspace*{\fill} \CommentCode{\textsf{NumPoissonGeometry} instance} \\
\LineCode \textsf{P = \{(1, 4): 1, (2, 5): 1, (3, 6): 1, (5, 6): `x2**2'\}}
    \hspace*{\fill} \CommentCode{dictionary for $\Pi$ in (\ref{EcTwist}) according to (\ref{EcMultivectorDic})} \\
\LineCode \textsf{alpha = \{(5,): 1\}}
    \hspace*{\fill} \CommentCode{dictionary for $\dd{x_{5}}$ according to (\ref{EcMultivectorDic})} \\
\LineCode \textsf{beta = \{(6,): 1\}}
    \hspace*{\fill} \CommentCode{dictionary for $\dd{x_{6}}$ according to (\ref{EcMultivectorDic})} \\
\LineCode \textsf{npg6.num\_one\_forms\_bracket(P, alpha, beta, mesh, pt\_output=True)} \newline
    \hspace*{\fill} \CommentCode{run \textsf{num\_one\_forms\_bracket} function with \textsf{pt\_output} flag}
    \tcblower
\textsf{tensor([\,\parbox[t]{\linewidth}{[[0.], [2.], [0.], [0.], [0.], [0.]], \ [[0.], [2.], [0.], [0.], [0.], [0.]], \newline
        $\ldots$, \newline
    [[0.], [2.], [0.], [0.], [0.], [0.]], \ [[0.], [2.], [0.], [0.], [0.], [0.]], \newline
    [[0.], [2.], [0.], [0.], [0.], [0.]], \ [[0.], [2.], [0.], [0.], [0.], [0.]]\,], \newline
    dtype=torch.float64)}}
\end{tcolorbox}
\end{example}

    \subsection{Gauge Transformations} \label{subsec:gauge}

Gauge transformations are used to simplify dynamical equations, they aid in reduction methods for dynamical systems. They are also used for the hamiltonization of nonholonomic systems \cite{GaugeNaranjo}. Recall that the {\em hamiltonization problem} consists in determining conditions that represent a dynamical system in Hamiltonian form.

Given a differential 2--form $\lambda$ on $M$, if the morphism \ec{\mathrm{id}_{\T^{\ast}{M}} - \lambda^{\flat} \circ \Pi^{\sharp}} is invertible, the $\lambda$--gauge transformation of a bivector field $\Pi$ on $M$ is the bivector field $\overline{\Pi}$ determined by the vector bundle morphism \ec{{\overline{\Pi}}^{\sharp} = \Pi^{\sharp} \circ (\mathrm{id}_{\T^{\ast}{M}} - \lambda^{\flat} \circ \Pi^{\sharp})^{-1}} \cite{SeveraTwist}. This morphism can be evaluated at a point in $M$ as follows:
    \begin{equation*}
        \overline{\Pi}^{\sharp}_{p} = -\big[\Pi^{ij}(p)\big]
        \left(\mathrm{I} - \big[\lambda_{kl}(p)\big]\big[\Pi^{ij}(p)\big]\right)^{-1}, \quad p \in M
     \end{equation*}
Here, \ec{[\Pi^{ij}]} is the matrix of $\Pi$ in (\ref{EcPiMatrix}), \ec{[\lambda_{kl}]} is the matrix of \ec{\lambda = 1/2\,\lambda_{kl}\dd{x^{k}} \wedge \dd{x^{l}}} and $\mathrm{I}$ denotes the \ec{m \times m} identity matrix, for \ec{k,l = 1,\ldots,m}. The morphism $\lambda^{\flat}:\T{M} \rightarrow \T^{\ast}{M}$ above is given by \ec{X \mapsto \ii_{X}\lambda}.

The function \textsf{num\_gauge\_transformation} evaluates the gauge transformation of a bivector field on a mesh in \ec{\RR{m}}.
\begin{algorithm}[H]
    \captionsetup{justification=centering}
    \caption{\ \textsf{num\_gauge\_transformation}(\emph{bivector, two\_form, mesh})} \label{AlgNumGauge}
        \rule{\textwidth}{0.4pt}
    \Input{a (Poisson) bivector field $\Pi$, a differential 2--form $\lambda$ and a mesh}
    \Output{evaluation of the gauge transformation of $\Pi$ induced by $\lambda$ at each point of the mesh \vspace{0.25cm}}
        \rule{\textwidth}{0.4pt}
    \begin{algorithmic}[1] 
        \Procedure{}{}
            \State $m$ $\gets$ dimension of the manifold
            \State \bluecolor{bivector} $\gets$ a variable encoding the (Poisson) bivector field $\Pi$
            \State \bluecolor{two\_form} $\gets$ a variable encoding the (Poisson) the differential 2--form $\lambda$
            \State \bluecolor{mesh} $\gets$ a $(k,m)$ array encoding the mesh
                \CommentNew{$k$: number of points in the mesh}
            \State \bluecolor{variable\_1} $\gets$ \textsf{num\_bivector\_to\_matrix}(\bluecolor{bivector}, \bluecolor{mesh})
                \CommentNew{see Algorithm \ref{AlgNumMatrixBivector}}
            \State \bluecolor{variable\_2} $\gets$ \textsf{num\_bivector\_to\_matrix}(\bluecolor{two\_form}, \bluecolor{mesh})
                \CommentNew{see Algorithm \ref{AlgNumMatrixBivector}}
            \State \bluecolor{variable\_3} $\gets$ the identity \ec{m \times m} matrix
            \State \bluecolor{variable\_4} $\gets$ an empty container
            \For {$0 \leq i < k$}
                \State \bluecolor{variable\_4}[$i$] $\gets$ \bluecolor{variable\_3} $-$ \bluecolor{variable\_2}[$i$] $*$ \bluecolor{variable\_1}[$i$]
            \Statex \CommentNew{matrix sum/multiplication}
            \EndFor
            \State \bluecolor{variable\_5} $\gets$ an empty container
            \For{$0 \leq i < k$}
                \If{\textsf{determinant}(\bluecolor{variable\_4}[$i$]) $\neq 0$}
                    \State \bluecolor{variable\_5}[$i$] $\gets$ \bluecolor{variable\_1}[$i$] $*$ \textsf{inverse}(\bluecolor{variable\_4}[$i$])
                    \CommentNew{matrix multiplication}
                \Else
                    \State \bluecolor{variable\_5}[$i$] $\gets$ \textsf{False}
                \EndIf
            \EndFor
            \State \rreturn{\bluecolor{variable\_5}}
        \EndProcedure
    \end{algorithmic}
\end{algorithm}
    \newpage

\begin{example}
For this example we will use the following result that we proved in \cite{CompuPoisson}:
\begin{proposition}\label{Prop:gauge}
Let $\Pi$ be a bivector field on a 3--dimensional smooth manifold $M$. Then, given a differential 2--form $\lambda$ on $M$, the $\lambda$--gauge transformation $\overline{\Pi}$ of $\Pi$ is well defined on the open subset,
    \begin{equation}\label{EcFdetGauge}
        \big\{F := 1 + \big\langle \lambda,\Pi \big\rangle \neq 0 \big\} \subseteq M.
    \end{equation}
Moreover, $\overline{\Pi}$ is given by
    \begin{equation*}\label{EcPiGauge3}
        \overline{\Pi} = \tfrac{1}{F}\Pi.
    \end{equation*}
In consequence, if $\Pi$ is Poisson, then $\overline{\Pi}$ is also Poisson.
\end{proposition}

Observe that if \ec{\langle \lambda,\Pi \rangle = 0} in (\ref{EcFdetGauge}), then $\Pi$ remains unchanged under the gauge transformation induced by $\lambda$. This holds for the Poisson bivector field \ec{\Pi_{\mathfrak{so}(3)}} in (\ref{EcPiSO3}) \cite{Newton} and the differential 2--form on \ec{\RR{3}} given by
    \begin{equation}\label{EcLambGague}
        \lambda = (x_1 - x_2)\dd{x_1}\wedge \dd{x_2} + (x_1 - x_3)\dd{x_1}\wedge \dd{x_3} + (x_2 - x_3)\dd{x_2}\wedge \dd{x_3}.
    \end{equation}
Then, for \ec{Q^{3}} being as in (\ref{EcCorners}), we can check the invariance of \ec{\Pi_{\mathfrak{so}(3)}} under $\lambda$ as follows:
\begin{tcolorbox}[arc=0mm, boxsep=0mm, skin=bicolor, colback=pink!15, colframe=blue!20, colbacklower=blue!0, breakable, halign=left]
\LineCode \textsf{npg3 = NumPoissonGeometry(3)}
    \hspace*{\fill} \CommentCode{\textsf{NumPoissonGeometry} instance} \\
\LineCode \textsf{P\_so3 = \big\{(1, 2): `x3', (1, 3): `-x2', (2, 3): `x1'\big\}}
    \hspace*{\fill} \CommentCode{dictionary for $\Pi_{\mathfrak{so}(3)}$ in (\ref{EcPiSO3}) according to (\ref{EcMultivectorDic})} \\
\LineCode \textsf{lambda = \big\{(1, 2): `x1 - x2', (1, 3): `x1 - x3', (2, 3): `x2 - x3'\big\}} \newline
    \hspace*{\fill} \CommentCode{dictionary for $\lambda$ in (\ref{EcLambGague}) according to (\ref{EcMultivectorDic})} \\
\LineCode \textsf{npg3.num\_gauge\_transformation(P\_so3, lambda, Qmesh, pt\_output=True)}
    \hspace*{\fill} \CommentCode{run \textsf{num\_gague\_transformation} function with \textsf{pt\_output} flag}
    \tcblower
\resizebox{\textwidth}{!}{
\textsf{tensor([\,\parbox[t]{\linewidth}{[[0., 0., \phantom{-}0.], [0., 0., 0.], [0., \phantom{-}0., 0.]], \
        [[0., 1., \phantom{-}0.], [-1., 0., 0.], [0., \phantom{-}0., 0.]], \newline
        [[0., 0., -1.], [0., 0., 0.], [1., \phantom{-}0., 0.]], \
        [[0., 1., -1.], [-1., 0., 0.], [1., \phantom{-}0., 0.]], \newline
        [[0., 0., \phantom{-}0.], [0., 0., 1.], [0., -1., 0.]], \
        [[0., 1., \phantom{-}0.], [-1., 0., 1.], [0., -1., 0.]], \newline
        [[0., 0., -1.], [0., 0., 1.], [1., -1., 0.]], \
        [[0., 1., -1.], [-1., 0., 1.], [1., -1., 0.]]\,], \newline
        dtype=torch.float64)}}
    }
\end{tcolorbox}
Notice that the output is the same as the second one in Example \ref{example:so3}, which encodes the evaluation of the bivector field $\Pi_{\mathfrak{so}(3)}$ at points of \ec{Q^{3}}.
\end{example}

    \subsection{Lie--Poisson Normal Forms on \ec{\RR{3}}} \label{subsec:normal}

Two Poisson bivector fields $\Pi$ and $\widetilde{\Pi}$ on $M$ are said to be equivalent (or isomorphic) if there exists a diffeomorphism \ec{F:M \rightarrow M} such that \ec{\widetilde{\Pi} = F^{\ast}\Pi}. Under this equivalence relation there exist 9 non--trivial normal forms of Lie--Poisson bivector fields on \ec{\RR{3}} \cite{LiuXU-92, Sheng}.

The function \textsf{num\_linear\_normal\_form\_R3} evaluates a normal form of a given Lie--Poisson bivector field on a mesh in \ec{\RR{3}}.

\begin{algorithm}[H]
    \captionsetup{justification=centering}
    \caption{\ \textsf{num\_linear\_normal\_form\_R3}(\emph{linear\_bivector, mesh})} \label{AlgNumNormal}
        \rule{\textwidth}{0.4pt}
    \Input{a Lie--Poisson bivector field $\Pi$ on $\RR{3}$ and a mesh}
    \Output{evaluation of a normal form of $\Pi$ at each point of the mesh}
        \rule{\textwidth}{0.4pt}
    \begin{algorithmic}[1] 
        \Procedure{}{}
            \State $m$ $\gets$ dimension of the manifold
            \State \bluecolor{linear\_bivector} $\gets$ a variable that encodes the Lie--Poisson bivector field $\Pi$
            \State \bluecolor{mesh} $\gets$ a $(k,m)$ array encoding the mesh
                \CommentNew{$k$: number of points in the mesh}
            \State \bluecolor{variable\_1} $\gets$ a variable encoding a linear normal form of \bluecolor{linear\_bivector}
            \Statex \CommentNew{We use the \textsf{linear\_normal\_form\_R3} method of \textsf{PoissonGeometry}}
            \State \rreturn{\textsf{num\_bivector}(\bluecolor{variable\_1}, \bluecolor{mesh})}
                \CommentNew{see Algorithm \ref{AlgNumBivector}}
        \EndProcedure
    \end{algorithmic}
\end{algorithm}

\begin{example}
Using the function \textsf{linear\_normal\_form\_R3} of \textsf{PoissonGeometry} we can verify that the Lie--Poisson bivector field on \ec{\RR{3}_{x}}
    \begin{equation}\label{EcPiNormal}
        \Pi =
        2(x_{2} + x_{3})\frac{\partial}{\partial x_1}\wedge \frac{\partial}{\partial x_2} +
        (x_{1} - x_{2})\frac{\partial}{\partial x_1}\wedge \frac{\partial}{\partial x_3} +
        (x_{1} + x_{2} + 2x_{3})\frac{\partial}{\partial x_2}\wedge \frac{\partial}{\partial x_3},
    \end{equation}
admits the following Poisson bivector field as a normal form:
    \begin{equation*}\label{EcPiNormal2}
        \Pi_{N} =
        (x_{1} - 4ax_{2})\frac{\partial}{\partial x_1}\wedge \frac{\partial}{\partial x_3} +
        (4ax_{1} + x_{2})\frac{\partial}{\partial x_2}\wedge \frac{\partial}{\partial x_3}, \quad a>0
    \end{equation*}
To evaluate this normal form of $\Pi$ at points of \ec{Q^{3}} (\ref{EcCorners}) we compute:
\begin{tcolorbox}[arc=0mm, boxsep=0mm, skin=bicolor, colback=pink!15, colframe=blue!20, colbacklower=blue!0, breakable, halign=left]
\LineCode \textsf{npg3 = NumPoissonGeometry(3)}
    \hspace*{\fill} \CommentCode{\textsf{NumPoissonGeometry} instance} \\
\LineCode \textsf{P = \big\{(1, 2): `2*(x2 + x3)', (1, 3): `x1 - x2', (2, 3): `x1 + x2 +2*x3'\big\}} \newline
    \hspace*{\fill} \CommentCode{dictionary for $\Pi$ in (\ref{EcPiNormal}) according to (\ref{EcMultivectorDic})} \\
\LineCode \textsf{npg3.num\_linear\_normal\_form\_R3(P, Qmesh, dict\_output=True)}
    \hspace*{\fill} \CommentCode{run \textsf{num\_linear\_normal\_form} function}
    \tcblower
\textsf{[\parbox[t]{\linewidth}{\{(1,3): \phantom{-}0.0, (2,3): 0.0\}, \ \{(1,3): \phantom{-}0.0, (2,3): 0.0\}, \
 \{(1,3): -4.0*a, (2,3): 1.0\}, \newline
    \{(1,3): -4.0*a, (2,3): 1.0\}, \{(1,3): 1.0, (2,3): 4.0*a\}, \{(1,3): 1.0, (2,3): 4.0*a\}, \newline
    \{(1,3): 1.0-4.0*a, (2,3): 4.0*a+1.0\}, \ \{(1,3): 1.0-4.0*a, (2,3): 4.0*a+1.0\}]}}
\end{tcolorbox}
\end{example}

The equivalence between $\Pi$ and \ec{\Pi_{N}} implies that the characteristic foliation of $\Pi$ is a open book foliation \cite{Obook}, as shown in Figure \ref{fig:OpenBook} below. In particular, $\Pi$ does not admit global, non--constant, Casimir functions.
\begin{figure}[H]
  \centering
  \includegraphics[width=0.25\textwidth]{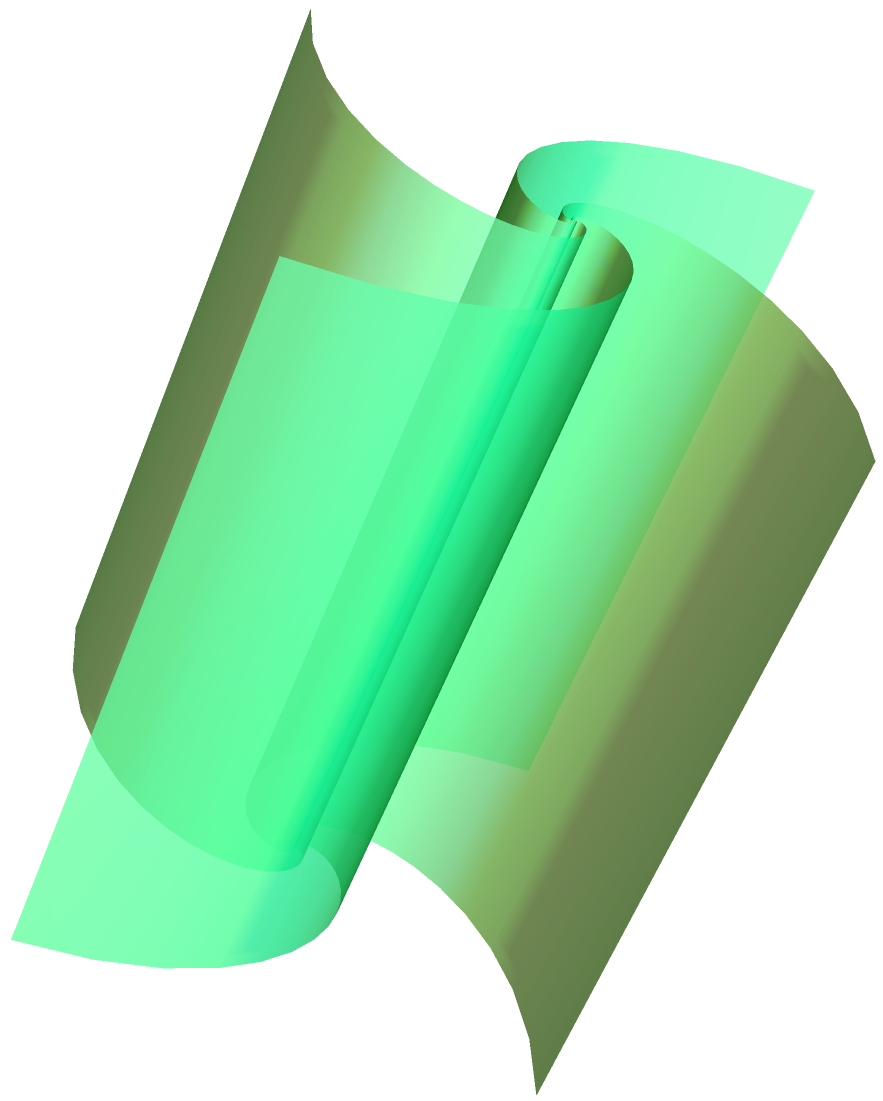}
  \caption{Symplectic (open book) foliation of $\Pi$ in (\ref{EcPiNormal}).} \label{fig:OpenBook}
\end{figure}

    \subsection{Flaschka--Ratiu Poisson Bivector Field} \label{subsec:fratiu}

If $M$ is a oriented manifold with volume form $\Omega$, the Poisson bivector field $\Pi$ determined by \ec{m-2} prescribed Casimir functions \ec{K_{1},...,K_{m-2}} on $M$, and defined by the formula
    \begin{equation*}
        \ii_{\Pi}\Omega = \dd{K_{1}} \wedge \cdots \wedge \dd{K_{m-2}},
    \end{equation*}
is called Flaschka--Ratiu bivector field \cite{GrabowskiFR, Damianou}. Observe that $\Pi$ is non--trivial on the open subset of $M$ where \ec{K_{1},\ldots,K_{m-2}} are (functionally) independent.

The function \textsf{num\_flaschka\_ratiu\_bivector} evaluates a Flaschka--Ratiu bivector field on a mesh in \ec{\RR{m}}.

\begin{algorithm}[H]
    \captionsetup{justification=centering}
    \caption{\ \textsf{num\_flaschka\_ratiu\_bivector}(\emph{casimirs\_list, mesh})} \label{AlgNumFRatiu}
        \rule{\textwidth}{0.4pt}
    \Input{a set of scalar functions and a mesh}
    \Output{evaluation of the Flaschka--Ratiu bivector field induced by these functions at each point of the mesh \vspace{0.25cm}}
        \rule{\textwidth}{0.4pt}
    \begin{algorithmic}[1] 
        \Procedure{}{}
            \State $m$ $\gets$ dimension of the manifold
            \State \bluecolor{casimirs\_list} $\gets$ a container with $m-2$ (string) expressions representing the set of scalar functions
                \CommentNew{each string expression represents a scalar function}
            \State \bluecolor{mesh} $\gets$ a $(k,m)$ array encoding the mesh
                \CommentNew{$k$: number of points in the mesh}
            \State \bluecolor{variable\_1} $\gets$ a variable encoding the Flaschka--Ratiu bivector field induced by \bluecolor{casimirs\_list}
                \CommentNew{We use the \textsf{flaschka\_ratiu\_bivector} method of \textsf{PoissonGeometry}}
            \State \rreturn{\textsf{num\_bivector}(\bluecolor{variable\_1}, \bluecolor{mesh})}
                \CommentNew{see Algorithm \ref{AlgNumBivector}}
        \EndProcedure
    \end{algorithmic}
\end{algorithm}

\begin{example}\label{example:so3FR}
Consider the following (Poisson) Flaschka--Ratiu bivector field on \ec{\RR{4}_{x}}
    \begin{equation*}
        \Pi =
        x_{3}\frac{\partial}{\partial x_1}\wedge \frac{\partial}{\partial x_2} -
        x_{2}\frac{\partial}{\partial x_1}\wedge \frac{\partial}{\partial x_3} -
        x_{1}\frac{\partial}{\partial x_2}\wedge \frac{\partial}{\partial x_3},
    \end{equation*}
that appears as a local model around singularities of a broken Lefschetz fibration on smooth 4--manifolds \cite{PabLef} (see, also \cite{PabloWrinFib, PabBott}). It is induced by the functions
    \begin{equation}\label{EcK1K2}
        K_{1} = \tfrac{1}{2}x_{4}, \quad K_{2} = - x_{1}^{2} + x_{2}^{2} + x_{3}^{2}.
    \end{equation}
To evaluate $\Pi$ at points of \ec{Q^{4}} (\ref{EcCorners}) we compute:
\begin{tcolorbox}[arc=0mm, boxsep=0mm, skin=bicolor, colback=pink!15, colframe=blue!20, colbacklower=blue!0, breakable, halign=left]
\LineCode \textsf{npg4 = NumPoissonGeometry(4)}
    \hspace*{\fill} \CommentCode{\textsf{NumPoissonGeometry} instance} \\
\LineCode \textsf{functions = [`1/2*x4', `-x1**2 + x2**2 + x3**2']} \newline
    \hspace*{\fill} \CommentCode{list containing string expressions for $K_{1}$ and $K_{2}$ in (\ref{EcK1K2}), in that order} \\
\LineCode \textsf{npg4.num\_flaschka\_ratiu\_bivector(functions, Qmesh, dict\_output=True)}
    \hspace*{\fill} \CommentCode{run \textsf{num\_flaschka\_ratiu\_bivector} function}
    \tcblower
\resizebox{\textwidth}{!}{
\textsf{[\parbox[t]{\linewidth}{\{(1, 2): 0.0, (1, 3): \phantom{-}0.0, (2, 3): \phantom{-}0.0\}, \
 \{(1, 2): 0.0, (1, 3): \phantom{-}0.0, (2, 3): \phantom{-}0.0\}, \newline
 \{(1, 2): 1.0, (1, 3): \phantom{-}0.0, (2, 3): \phantom{-}0.0\}, \
 \{(1, 2): 1.0, (1, 3): \phantom{-}0.0, (2, 3): \phantom{-}0.0\}, \newline
 \{(1, 2): 0.0, (1, 3): -1.0, (2, 3): \phantom{-}0.0\}, \
 \{(1, 2): 0.0, (1, 3): -1.0, (2, 3): \phantom{-}0.0\}, \newline
 \{(1, 2): 1.0, (1, 3): -1.0, (2, 3): \phantom{-}0.0\}, \
 \{(1, 2): 1.0, (1, 3): -1.0, (2, 3): \phantom{-}0.0\}, \newline
 \{(1, 2): 0.0, (1, 3): \phantom{-}0.0, (2, 3): -1.0\}, \
 \{(1, 2): 0.0, (1, 3): \phantom{-}0.0, (2, 3): -1.0\}, \newline
 \{(1, 2): 1.0, (1, 3): \phantom{-}0.0, (2, 3): -1.0\}, \
 \{(1, 2): 1.0, (1, 3): \phantom{-}0.0, (2, 3): -1.0\}, \newline
 \{(1, 2): 0.0, (1, 3): -1.0, (2, 3): -1.0\}, \
 \{(1, 2): 0.0, (1, 3): -1.0, (2, 3): -1.0\}, \newline
 \{(1, 2): 1.0, (1, 3): -1.0, (2, 3): -1.0\}, \
 \{(1, 2): 1.0, (1, 3): -1.0, (2, 3): -1.0\}]}}
    }
\end{tcolorbox}
\end{example}

    \section{Algorithmic Complexity and Performance} \label{sec:ComplexityPerformance}

In this section we will present an \emph{approximation} to the (worst--case) time complexity of the twelve algorithms in Table \ref{table:Funs-Algos-Exes}, as well as a time performance analysis of our Python implementation of the functions in the module \textsf{NumPoissonGeometry}. Recall that we are assuming the input functions to be computable.

    \subsection{Complexity}

The time complexity of certain algorithms depends on the nature and structure of the input data. In our methods two important processes in the implementation depend on the \emph{length} of the input data items: converting string expressions to symbolic variables, and transforming these symbolic expressions into functions that allow a numerical evaluation. Therefore, for the analysis of our algorithms we define,
    \begin{equation*}
        |C| := \underset{x \,\in\, C}{\operatorname{max}}\,\{\mathrm{len}(x)\},
    \end{equation*}
for a container $C$ with (string expression) items that encode the coefficients (scalar functions) of the coordinate expression of a multivector field or a differential form. For example, as illustrated in (\ref{EcMultivectorDic}), we use dictionaries for our Python implementation. In this case, $C$ is a tuple (or list) that contains all the values of a such dictionary.

In the following Table \ref{table:complexity} we record the approximate worst--case time complexities for the twelve methods of our \textsf{NumPoissonGeometry} module:

\begin{table}[H]
    \centering
    \caption{Worst--case time complexity of \textsf{NumPoissonGeometry} methods. In the second column: $m$ denotes the dimension of \ec{\RR{m}}, $k$ is the number of points in a mesh on \ec{\RR{m}}, we denote by \ec{[\cdot]} the integer part function and by $\mathrm{comb}$ a combination.} \label{table:complexity}
    \resizebox{\textwidth}{!}{
 \begin{tabular}{|l|l|}
    \hline
        \multicolumn{1}{|c|}{\textbf{Method}} & \multicolumn{1}{c|}{\textbf{Time Complexity}} \\
    \hline
    \hline
        \hyperref[AlgNumBivector]{\phantom{0}1. \textsf{num\_bivector\_field}}
        & \ec{\OO(m^2k|bivector|)} \\
    \hline
        \hyperref[AlgNumMatrixBivector]{\phantom{0}2. \textsf{num\_bivector\_to\_matrix}}
        & \ec{\OO(m^{2}k|bivector|)} \\
    \hline
        \hyperref[AlgNumHamVF]{\phantom{0}3. \textsf{num\_hamiltonian\_vf}}
        & \ec{\OO(mk(m|bivector| + \mathrm{len}(ham\_function)))} \\
    \hline
        \hyperref[AlgNumPoissonBracket]{\phantom{0}4. \textsf{num\_poisson\_bracket}}
        & \ec{\OO(mk(m|bivector| + \mathrm{len}(function\_1) + \mathrm{len}(function\_2)))} \\
    \hline
        \hyperref[AlgNumSharp]{\phantom{0}5. \textsf{num\_sharp\_morphism}}
        & \ec{\OO(m(mk|bivector| + |one\_form|))} \\
    \hline
        \hyperref[AlgNumCoboundary]{\phantom{0}6. \textsf{num\_coboundary\_operator}}
        & \ec{\OO(\mathrm{comb}(m,[m/2])|bivector|\mathrm{len}(function)(m^{5} + k))} \\
    \hline
        \hyperref[AlgNumModularVF]{\phantom{0}7. \textsf{num\_modular\_vf}}
        & \ec{\OO(\mathrm{comb}(m,[m/2])|bivector|\mathrm{len}(function)(m+k))} \\
    \hline
        \hyperref[AlgNumCurl]{\phantom{0}8. \textsf{num\_curl\_operator}}
        & \ec{\OO(\mathrm{comb}(m,[m/2])|multivector|\mathrm{len}(function)(m+k))} \\
    \hline
        \hyperref[AlgNumOneFormsB]{\phantom{0}9. \textsf{num\_one\_forms\_bracket}}
        & \ec{\OO(m^2k|bivector||one\_form\_1||one\_form\_2|)} \\
    \hline
        \hyperref[AlgNumGauge]{10. \textsf{num\_gauge\_transformation}}
        & \ec{\OO\big(m^2k(m^5 + |bivector| + |two\_form|)\big)} \\
    \hline
        \hyperref[AlgNumNormal]{11. \textsf{num\_linear\_normal\_form\_R3}}
        & \ec{\OO(k|bivector|)} \\
    \hline
        \hyperref[AlgNumFRatiu]{12. \textsf{num\_flaschka\_ratiu\_bivector}}
        & \ec{\OO(m^{6}k|bivector|)} \\
    \hline
    \hline
 \end{tabular}
    }
\end{table}

\begin{remark}
Observe that the time complexities in Table \ref{table:complexity} depend linearly on the number of points in the mesh ($k$).
\end{remark}

Recall that the time complexity of independent processes equals the sum of the respective time complexities of each process. With this in mind, we present the deduction of the time complexities of our numerical methods presented in Table \ref{table:complexity}.

    \subsubsection{Polynomial Complexity}

The following methods have polynomial time complexities.

\begin{lemma}
The time complexity of the \hyperref[AlgNumBivector]{\textsf{num\_bivector\_field}} method is approximately
    \begin{equation*}
        \OO(m^2k|bivector|).
    \end{equation*}
\end{lemma}
\begin{proof}
Consider the Algorithm \ref{AlgNumBivector} with \emph{bivector} and \emph{mesh} inputs. The time complexity of our implementation depends on:
    \begin{enumerate}[label=\roman*.]
       \item \textsf{Line 5} $\rightarrow$ \ec{\OO(m^{2} |bivector|)}: we iterate over \emph{bivector}, and the transformation of each \emph{bivector} item depends on its length.

       \item \textsf{Line 6} $\rightarrow$ \ec{\OO(km^{2} |bivector|)}: we iterate over \emph{mesh}, and the evaluation of \emph{bivector} depends on its length and on the length of its items.
    \end{enumerate}
Hence, the time complexity of \textsf{num\_bivector\_field} is \ec{\OO(m^2k|bivector|)}.
\end{proof}


\begin{lemma}
The time complexity of the \hyperref[AlgNumMatrixBivector]{\textsf{num\_bivector\_to\_matrix}} method is approximately
    \begin{equation*}
        \OO(m^{2}k|bivector|).
    \end{equation*}
\end{lemma}
\begin{proof}
Consider the Algorithm \ref{AlgNumMatrixBivector} with \emph{bivector} and \emph{mesh} inputs. The time complexity of our implementation depends on:
    \begin{enumerate}[label=\roman*.]
       \item \textsf{Line 5} $\rightarrow$ \ec{\OO(m^{2}|bivector|)}

       \item \textsf{Line 6} $\rightarrow$ \ec{\OO(m^{2}|bivector|)}: we iterate over the \ec{m \times m} matrix \emph{variable\_1}, and the transformation of each \emph{bivector} item depends on its length.

       \item \textsf{Line 7} $\rightarrow$ \ec{\OO(km^{2}|bivector|)}: we iterate over \emph{mesh}, and the evaluation of \emph{variable\_1} depends on its dimension and on the length of the \emph{bivector} items.
    \end{enumerate}
Hence, the time complexity of \textsf{num\_bivector\_to\_matrix} is \ec{\OO(m^{2}k|bivector|)}.
\end{proof}


\begin{lemma}
The time complexity of the \hyperref[AlgNumHamVF]{\textsf{num\_hamiltonian\_vf}} method is approximately
    \begin{equation*}
        \OO(mk(m|bivector| + \mathrm{len}(ham\_function)))
    \end{equation*}
\end{lemma}
\begin{proof}
Consider the Algorithm \ref{AlgNumHamVF} with \emph{bivector}, \emph{ham\_function} and \emph{mesh} inputs. The time complexity of our implementation depends on:
    \begin{enumerate}[label=\roman*.]

       \item \textsf{Line 6} $\rightarrow$ \ec{\OO(m\mathrm{len}(ham\_function))}: calculate the gradient of \emph{ham\_function} depends on the length of this expression and on the dimension of \ec{\RR{m}}.

       \item \textsf{Line 7} $\rightarrow$ \ec{\OO(m\mathrm{len}(ham\_function))}: we iterate over the \ec{m \times 1} matrix \emph{variable\_1}, and the transformation of each of its items depends on the length of \emph{ham\_function}.

       \item \textsf{Line 8} $\rightarrow$ \ec{\OO(km\mathrm{len}(ham\_function))}: we iterate over \emph{mesh}, and the evaluation of \emph{variable\_1} depends on its dimension and on the length of \emph{ham\_function}.

       \item \textsf{Line 9} $\rightarrow$ \ec{\OO(m^{2}k|bivector|)}: the \hyperref[AlgNumMatrixBivector]{\textsf{num\_bivector\_to\_matrix}} method has time complexity \ec{\OO(m^{2}k|bivector|)}.

       \item \textsf{Lines 11-13} $\rightarrow$ \ec{\OO(km^{2})}: we iterate over the set of indices \ec{\{0, \ldots, k-1\}}, and the product of the \emph{variable\_2} and \emph{variable\_3} items has time complexity \ec{\OO(m^{2})}.
    \end{enumerate}
Hence, the time complexity of \textsf{num\_hamiltonian\_vf} is
    \begin{equation*}
        \OO(m k (m |bivector| + \mathrm{len} (ham\_function))).
    \end{equation*}
\end{proof}


\begin{lemma}
The time complexity of the \hyperref[AlgNumPoissonBracket]{\textsf{num\_poisson\_bracket}} method is approximately
    \begin{equation*}
        \OO(mk(m|bivector| + \mathrm{len}(function\_1) + \mathrm{len}(function\_2))).
    \end{equation*}
\end{lemma}
\begin{proof}
Consider the Algorithm \ref{AlgNumPoissonBracket} with \emph{bivector}, \emph{function\_1}, \emph{function\_2} and \emph{mesh} inputs. The time complexity of our implementation depends on:
    \begin{enumerate}[label=\roman*.]

       \item \textsf{Line 9} $\rightarrow$ \ec{\OO(m\mathrm{len}(function\_2))}: calculate the gradient of \emph{function\_2} depends on the length of this expression and on the dimension of \ec{\RR{m}}.

       \item \textsf{Line 10} $\rightarrow$ \ec{\OO(m\mathrm{len}(function\_2))}: we iterate over the \ec{m \times 1} matrix \emph{variable\_1}, and the transformation of each of its items depends on the length of \emph{function\_2}.

       \item \textsf{Line 11} $\rightarrow$ \ec{\OO(km\mathrm{len}(function\_2))}: we iterate over \emph{mesh}, and the evaluation of \emph{variable\_1} depends on its dimension and on the length of \emph{function\_2}.

       \item \textsf{Line 12} $\rightarrow$ \ec{\OO(mk(m|bivector| + \mathrm{len}(function\_1)))}: the \hyperref[AlgNumHamVF]{\textsf{num\_hamiltonian\_vf}} has time complexity \ec{\OO(mk(m|bivector| + \mathrm{len}(ham\_function)))}.

       \item \textsf{Lines 14-16} $\rightarrow$ \ec{\OO(km^{2})}: we iterate over the set of indices \ec{\{0, \ldots, k-1\}}, and the product of the \emph{variable\_2} and \emph{variable\_3} items has time complexity \ec{\OO(m^{2})}.
    \end{enumerate}
Hence, the time complexity of \textsf{num\_poisson\_bracket} is
\[
\ec{\OO(mk(m|bivector| + \mathrm{len}(function\_1) + \mathrm{len}(function\_2)))}.
\]
\end{proof}


\begin{lemma}
The time complexity of the \hyperref[AlgNumSharp]{\textsf{num\_sharp\_morphism}} method is approximately
    \begin{equation*}
        \OO(m(mk|bivector| + |one\_form|)).
    \end{equation*}
\end{lemma}
\begin{proof}
Consider the Algorithm \ref{AlgNumSharp} with \emph{bivector}, \emph{one\_form} and \emph{mesh} inputs. The time complexity of our implementation depends on:
    \begin{enumerate}[label=\roman*.]
       \item \textsf{Line 7} $\rightarrow$ \ec{\OO(m|one\_form|)}: we iterate over the \ec{(m,1)} container \emph{variable\_1}, and the transformation of each of its items depends on the length of \emph{one\_form} items.

       \item \textsf{Line 9} $\rightarrow$ \ec{\OO(m^{2}k|bivector|)}: the \hyperref[AlgNumMatrixBivector]{\textsf{num\_bivector\_to\_matrix}} method has time complexity \ec{\OO(m^{2}k|bivector|)}.

       \item \textsf{Lines 11-13} $\rightarrow$ \ec{\OO(km^{2})}: we iterate over the set of indices \ec{\{0, \ldots, k-1\}}, and the product of the \emph{variable\_2} and \emph{variable\_3} items has time complexity \ec{\OO(m^{2})}.
    \end{enumerate}
Hence, the time complexity of \textsf{num\_sharp\_morphism} is
    \begin{equation*}
        \OO(m (m k |bivector| + |one\_form|)).
    \end{equation*}
\end{proof}


\begin{lemma}
The time complexity of the \hyperref[AlgNumOneFormsB]{\textsf{num\_one\_forms\_bracket}} method is approximately
    \begin{equation*}
        \OO(m^2k|bivector| \\ |one\_form\_1||one\_form\_2|)
    \end{equation*}
\end{lemma}
\begin{proof}
Consider the Algorithm \ref{AlgNumOneFormsB} with \emph{bivector}, \emph{one\_form\_1}, \emph{one\_form\_2} and \emph{mesh} inputs. The time complexity of our implementation depends on:
    \begin{enumerate}[label=\roman*.]

       \item \textsf{Line 10} $\rightarrow$ \ec{\OO(m^2(|one\_form\_1| + |one\_form\_2|))}: calculate the Jacobian matrix of \emph{variable\_i} depends on the length of \emph{one\_form\_1} and \emph{one\_form\_2} items and on the dimension of \ec{\RR{m}}.

       \item \textsf{Line 11} $\rightarrow$ \ec{\OO(m^2(|one\_form\_1| + |one\_form\_2|))}: we iterate over the \ec{m \times m} matrix \emph{variable\_3\_i}, and the transformation of each of its items depends on the length of \emph{one\_form\_1} and \emph{one\_form\_2} items.

       \item \textsf{Line 13} $\rightarrow$ \ec{\OO(m(mk|bivector| + |one\_form\_1| + |one\_form\_2|))}: the \hyperref[AlgNumSharp]{\textsf{num\_sharp\_ morphism}} method has time complexity \ec{\OO(m(mk|bivector| + |one\_form|))}.

       \item \textsf{Lines 16-19} $\rightarrow$ \ec{\OO(km^{2})}: we iterate over the set of indices \ec{\{0, \ldots, k-1\}}, and the product of the \emph{variable\_4\_i} and \emph{variable\_5\_j} items has time complexity \ec{\OO(m^{2})}, for \ec{i,j = 1,2}.

       \item \textsf{Line 20} $\rightarrow$ \ec{\OO(m^2|bivector||one\_form\_1|)}: the \textsf{sharp\_morphism} method has time complexity \ec{\OO(m^2 |bivector||one\_form|)}.

       \item \textsf{Line 22} $\rightarrow$ \ec{\OO(m|one\_form\_1||one\_form\_2|)}: calculate the gradient of \emph{variable\_9} depends on the length of \emph{one\_form\_1} and \emph{one\_form\_2} items, and on the dimension of \ec{\RR{m}}.

       \item \textsf{Line 23} $\rightarrow$ \ec{\OO(m(|one\_form\_1||one\_form\_2|))}: we iterate over the \ec{m \times 1} matrix \emph{variable\_10},  and the transformation of each of its items depends on the length of \emph{one\_form\_1} and \emph{one\_form\_2} items.

       \item \textsf{Lines 26-28} $\rightarrow$ \ec{\OO(km)}: we iterate over the set of indices \ec{\{0, \ldots, k-1\}}, and the sum of variables \emph{variable\_6}, \emph{variable\_7} and \emph{variable\_11} items has time complexity \ec{\OO(m)}.
    \end{enumerate}
Hence, the time complexity of \textsf{num\_one\_forms\_bracket} is
\[
\OO( m^2 k |bivector| |one\_form\_1| |one\_form\_2|).
\]
\end{proof}


\begin{lemma}
The time complexity of the \hyperref[AlgNumGauge]{\textsf{num\_gauge\_transformation}} method is approximately \ec{\OO(m^2k(m^5 + |bivector| + |two\_form|))}.
\end{lemma}
\begin{proof}
Consider the Algorithm \ref{AlgNumGauge} with \emph{bivector}, \emph{two\_form} and \emph{mesh} inputs. The time complexity of our implementation depends on:
    \begin{enumerate}[label=\roman*.]
       \item \textsf{Line 6} $\rightarrow$ \ec{\OO(m^{2}k|bivector|)}: the \hyperref[AlgNumMatrixBivector]{\textsf{num\_bivector\_to\_matrix}} method has time complexity \ec{\OO(m^{2}k|bivector|)}.

       \item \textsf{Line 7} $\rightarrow$ \ec{\OO(m^{2}k|two\_form|)}: the \hyperref[AlgNumMatrixBivector]{\textsf{num\_bivector\_to\_matrix}} method has complexity \newline \ec{\OO(m^{2}k|bivector|)}.

       \item \textsf{Lines 10-12} $\rightarrow$ \ec{\OO(km^{2})}: we iterate over the set of indices \ec{\{0, \ldots, k-1\}}, and the matrix operations between the \emph{variable\_1}, \emph{variable\_2} and \emph{variable\_3} items has time complexity \ec{\OO(m^{2})}.

       \item \textsf{Lines 14-20} $\rightarrow$ \ec{\OO(km^{10})}: we iterate over the set of indices \ec{\{0, \ldots, k-1\}}, calculate the determinant of each item of \emph{variable\_4} has time complexity \ec{\OO(m^4)}, the inverse has time complexity \ec{\OO(m^3)}, and the matrix product in \textsf{line 16} has time complexity \ec{\OO(m^3)}.
    \end{enumerate}
Hence, the time complexity of \textsf{num\_gauge\_transformation} is
    \begin{equation*}
        \OO(m^2 k (m^5 + |bivector| + |two\_form|)).
    \end{equation*}
\end{proof}


\begin{lemma}
The time complexity of the \hyperref[AlgNumNormal]{\textsf{num\_linear\_normal\_form\_R3}} method is approximately \ec{\OO(k|bivector|)}.
\end{lemma}
\begin{proof}
Consider the Algorithm \ref{AlgNumNormal} with \emph{linear\_bivector} and \emph{mesh} inputs. The time complexity of our implementation depends on:
    \begin{enumerate}[label=\roman*.]
       \item \textsf{Line 5} $\rightarrow$ \ec{\OO(|bivector|)}: the \textsf{linear\_normal\_form\_R3} method has time complexity \ec{\OO(|bivector|)}.

       \item \textsf{Line 6} $\rightarrow$ \ec{\OO(k|bivector|)}: the \hyperref[AlgNumBivector]{\textsf{num\_bivector\_field}} method has time complexity \ec{\OO(m^{2}k|bivector|)}.
    \end{enumerate}
Hence, the time complexity of \textsf{num\_linear\_normal\_form\_R3} is \ec{\OO(k|bivector|)}.
\end{proof}


\begin{lemma}
The time complexity of the \hyperref[AlgNumFRatiu]{\textsf{num\_flaschka\_ratiu\_bivector}} method is approximately \ec{\OO(m^{6}k|bivector|)}.
\end{lemma}
\begin{proof}
Consider the Algorithm \ref{AlgNumFRatiu} with \emph{casimir\_list} and \emph{mesh} inputs. The time complexity of our implementation depends on:
    \begin{enumerate}[label=\roman*.]
       \item \textsf{Line 5} $\rightarrow$ \ec{\OO(m^{6}|casimir\_list|)}: the \textsf{flaschka\_ratiu\_bivector} method has complexity \newline \ec{\OO(m^{6}|casimir\_list|)}.

       \item \textsf{Line 6} $\rightarrow$ \ec{\OO(m^{2}k|casimir\_list|)}: the \hyperref[AlgNumBivector]{\textsf{num\_bivector\_field}} method has time complexity \newline \ec{\OO(m^{2}k|bivector|)}.
    \end{enumerate}
Hence, the time complexity is \ec{\OO(m^{6}k|bivector|)}.
\end{proof}

    \subsubsection{Exponential Complexity}

Due to the nature of the concatenated loops in our methods, and because we need to calculate ordered index permutations, the following methods have exponential time complexities.

\begin{lemma}
The time complexity of the \hyperref[AlgNumCoboundary]{\textsf{num\_coboundary\_operator}} method is approximately \ec{\OO(\mathrm{comb}(m,[m/2]) |bivector|\mathrm{len}(function)(m^{5} + k))}.
\end{lemma}
\begin{proof}
Consider the Algorithm \ref{AlgNumCoboundary} with \emph{bivector}, \emph{multivector} and \emph{mesh} inputs. The time complexity of our implementation depends on:
    \begin{enumerate}[label=\roman*.]
       \item \textsf{Line 6} $\rightarrow$ \ec{\OO(m^{5}\mathrm{comb}(m,[m/2])|bivector||multivector|)}: the \textsf{lichnerowicz\_poisson\_ operator} method has time complexity
           \begin{equation*}
                \OO(m^{5} \mathrm{comb} (m, [m/2]) |bivector||multivector|).
           \end{equation*}

       \item \textsf{Line 7} $\rightarrow$ \ec{\OO(\mathrm{comb}(m,[m/2])|bivector||multivector|)}: we iterate over \emph{variable\_1}, and the transformation of each of its items depends on the length of the \emph{bivector} and \emph{multivector} items.

       \item \textsf{Line 8} $\rightarrow$ \ec{\OO(k\mathrm{comb}(m,[m/2])|bivector||multivector|)}: we iterate over \emph{mesh}, and the evaluation of \emph{variable\_1} depends on its length and on the length of the \emph{bivector} and \emph{multivector} items.
    \end{enumerate}
Hence, the time complexity of \textsf{num\_coboundary\_operator} is \[ \ec{\OO(\mathrm{comb}(m,[m/2])|bivector|\mathrm{len}(function)(m^{5} + k))}. \]
\end{proof}


\begin{lemma}
The time complexity of the \hyperref[AlgNumModularVF]{\textsf{num\_modular\_vf}} method is approximately
    \begin{equation*}
        \OO(\mathrm{comb}(m,[m/2]) |bivector|\mathrm{len}(function)(m+k)).
    \end{equation*}
\end{lemma}
\begin{proof}
Consider the Algorithm \ref{AlgNumModularVF} with \emph{bivector}, \emph{function} and \emph{mesh} inputs. The time complexity of our implementation depends on:
    \begin{enumerate}[label=\roman*.]
       \item \textsf{Line 6} $\rightarrow$ \ec{\OO(m\,\mathrm{comb}(m,[m/2])|bivector| \mathrm{len}(function))}: the \textsf{modular\_vf} method has time complexity \ec{\OO(m\,\mathrm{comb}(m,[m/2])|bivector| \mathrm{len}(function))}.

       \item \textsf{Line 7} $\rightarrow$ \ec{\OO(\mathrm{comb}(m,[m/2])|bivector|\mathrm{len}(function))}: we iterate over \emph{variable\_1}, and the transformation of each of its items depends on the length of the \emph{bivector} items and the \emph{function} string expression.

       \item \textsf{Line 8} $\rightarrow$ \ec{\OO(k\mathrm{comb}(m,[m/2]|bivector|\mathrm{len}(function)))}: we iterate over \emph{mesh}. The evaluation of \emph{variable\_1} depends on its length, and on the length of the \emph{bivector} items and the \emph{function} string expression.
    \end{enumerate}
Hence, the time complexity of \textsf{num\_modular\_vf} is
\[
\OO(\mathrm{comb} (m, [m/2]) |bivector| \mathrm{len}(function)(m + k)).
\]
\end{proof}


\begin{lemma}
The time complexity of the \hyperref[AlgNumCurl]{\textsf{num\_curl\_operator}} method is approximately
    \begin{equation*}
        \OO(\mathrm{comb}(m,[m/2]) |multivector|\mathrm{len}(function)(m+k)).
    \end{equation*}
\end{lemma}
\begin{proof}
Consider the Algorithm \ref{AlgNumCurl} with \emph{bivector}, \emph{function} and \emph{mesh} inputs. The time complexity of our implementation depends on:
    \begin{enumerate}[label=\roman*.]
       \item \textsf{Line 6} $\rightarrow$ \ec{\OO(m\,\mathrm{comb}(m,[m/2])|multivector| \mathrm{len}(function))}: the \textsf{curl\_operator} method has time complexity \ec{\OO(m\,\mathrm{comb}(m,[m/2])|multivector| \mathrm{len}(function))}.

       \item \textsf{Line 7} $\rightarrow$ \ec{\OO(\mathrm{comb}(m,[m/2])|multivector|\mathrm{len}(function))}: we iterate over \emph{variable\_1}, and the transformation of each of its items depends on the length of the \emph{multivector} items and the \emph{function} string expression.

       \item \textsf{Line 8} $\rightarrow$ \ec{\OO(k\mathrm{comb}(m,[m/2]|multivector|\mathrm{len}(function)))}: we iterate over \emph{mesh}. The evaluation of \emph{variable\_1} depends on its length, and on the length of the \emph{multivector} items and the \emph{function} string expression.
    \end{enumerate}
Hence, the time complexity of \textsf{num\_curl\_operator} is \[ \ec{\OO(\mathrm{comb}(m,[m/2])|multivector|\mathrm{len}(function)(m+k))}.\]
\end{proof}

    \subsection{Performance}

The time performance of each function in \textsf{NumPoissonGeometry} was experimentally measured by evaluating concrete examples on \ec{\RR{2}}, and \ec{\RR{3}}, on precalculated (irregular) meshes with \ec{10^{\kappa}} points, for \ec{\kappa = 3,\ldots,7}. These meshes were generated by means of random samples extracted from a uniform distribution in the interval \ec{[0,1)}.

All the numerical experiments were performed on a workstation equipped with 48 GB of main memory in a 3 $\times$ 16 GB ddr4 module configuration and an Intel(R) Core(TM) i7-6700 CPU @ 3.40GHz CPU, running at 3.4 GHz for a peak theoretical performance of four cores.

    \subsubsection{Two Dimensional Case}

For the performance tests of our methods in dimension two, we have set the non--degenerate Poisson bivector field
    \begin{equation*}
        \Pi_{0} = \frac{\partial}{\partial x_1}\wedge \frac{\partial}{\partial x_2},
    \end{equation*}
induced by the standard symplectic structure \ec{\dd{x_{1}} \wedge \dd{x_{2}}} on \ec{\RR{2}_{x}}.

\begin{table}[H]
    \centering
    \caption{Input data used for the time performance tests of functions 1-10 in Table \ref{table:performanceR2}.}
    \resizebox{\textwidth}{!}{
 \begin{tabular}{|c|l|c|l|}
        \hline
    \multicolumn{1}{|c|} {\textbf{Function}} & \multicolumn{1}{c|}{\textbf{Input}} & \multicolumn{1}{|c|} {\textbf{Function}} & \multicolumn{1}{c|}{\textbf{Input}} \\
        \hline
        \hline
    1 & \ec{\Pi_{0}} &
    6 & \ec{\Pi_{0}}, \ec{W = x_{2} \frac{\partial}{\partial x_{1}} - x_{1} \frac{\partial}{\partial x_{2}}} \\
        \hline
    2 & \ec{\Pi_{0}} &
    7 & \ec{\Pi_{0}}, \ec{f=1} \\
        \hline
    3 & \ec{\Pi_{0}}, \ec{h = x_{1}^{2} + x_{2}^{2}} &
    8 & \ec{\Pi_{0}}, \ec{f=1} \\
        \hline
    4 & \ec{\Pi_{0}}, \ec{f = x_{1}^{2} + x_{2}^{2}}, \ec{g=x_{1} + x_{2}} &
    9 & \ec{\Pi_{0}}, \ec{\alpha = x_{1}\dd{x_{1}} + x_{2}\dd{x_{2}}}, \ec{\beta = \dd{x_{1}} + \dd{x_{2}}} \\
        \hline
    5 & \ec{\Pi_{0}}, \ec{\alpha = x_{1}\dd{x_{1}} + x_{2}\dd{x_{2}}} &
    10 & \ec{\Pi_{0}}, \ec{\lambda = \dd{x_{1}} \wedge \dd{x_{2}}} \\
        \hline
        \hline
 \end{tabular}
    }
 \end{table}

Table \ref{table:performanceR2} lists the mean time in seconds (with standard deviation) it takes to evaluate the first ten functions in \textsf{NumPoissonGeometry} (see, Table \ref{table:Funs-Algos-Exes}) on a irregular mesh on \ec{\RR{3}_{x}} with \ec{10^{\kappa}} points, computed by taking twenty-five samples, for \ec{\kappa = 3,\ldots,7}.

\begin{table}[H]
    \centering
    \caption{Summary of the time performance of \textsf{NumPoissonGeometry} functions, in dimension two.} \label{table:performanceR2}
    \resizebox{\textwidth}{!}{
 \begin{tabular}{|l|l l l l l|}
    \hline
    \hline
        \multirow{2}{*}{\hspace{1.25cm}\textbf{Function}} & \multicolumn{5}{c|}{\textbf{Points in mesh/Processing time (in seconds)}} \\
        & \hspace{0.85cm}$10^{3}$ & \hspace{0.85cm}$10^{4}$ & \hspace{0.85cm}$10^{5}$ & \hspace{0.85cm}$10^{6}$ & \hspace{0.85cm}$10^{7}$ \\
    \hline
        \hyperref[AlgNumBivector]{\phantom{0}1. \textsf{num\_bivector\_field}}
        &   0.004 \textcolor{gray}{$\pm$ 0.689}
        &   0.038 \textcolor{gray}{$\pm$ 0.009}
        &   0.356 \textcolor{gray}{$\pm$ 0.002}
        &   \phantom{0}3.545 \textcolor{gray}{$\pm$ 0.026}
        &   \phantom{0}35.711 \textcolor{gray}{$\pm$ 0.164}  \\
    \hline
        \hyperref[AlgNumMatrixBivector]{\phantom{0}2. \textsf{num\_bivector\_to\_matrix}}
        &   0.006 \textcolor{gray}{$\pm$ 3.633}
        &   0.046 \textcolor{gray}{$\pm$ 0.001}
        &   0.438 \textcolor{gray}{$\pm$ 0.004}
        &   \phantom{0}4.442 \textcolor{gray}{$\pm$ 0.037}
        &   \phantom{0}45.155 \textcolor{gray}{$\pm$ 1.466} \\
    \hline
        \hyperref[AlgNumHamVF]{\phantom{0}3. \textsf{num\_hamiltonian\_vf}}
        &   0.014 \textcolor{gray}{$\pm$ 0.001}
        &   0.112 \textcolor{gray}{$\pm$ 0.006}
        &   1.096 \textcolor{gray}{$\pm$ 0.021}
        &   10.867 \textcolor{gray}{$\pm$ 0.044}
        &   108.460 \textcolor{gray}{$\pm$ 0.726} \\
    \hline
        \hyperref[AlgNumPoissonBracket]{\phantom{0}4. \textsf{num\_poisson\_bracket}}
        &   0.021 \textcolor{gray}{$\pm$ 0.006}
        &   0.169 \textcolor{gray}{$\pm$ 0.001}
        &   1.652 \textcolor{gray}{$\pm$ 0.008}
        &   16.721 \textcolor{gray}{$\pm$ 0.049}
        &   168.110 \textcolor{gray}{$\pm$ 1.637} \\
    \hline
        \hyperref[AlgNumSharp]{\phantom{0}5. \textsf{num\_sharp\_morphism}}
        &   0.014 \textcolor{gray}{$\pm$ 0.658}
        &   0.111 \textcolor{gray}{$\pm$ 0.001}
        &   1.068 \textcolor{gray}{$\pm$ 0.007}
        &   10.725 \textcolor{gray}{$\pm$ 0.142}
        &   107.275 \textcolor{gray}{$\pm$ 0.667} \\
    \hline
        \hyperref[AlgNumCoboundary]{\phantom{0}6. \textsf{num\_coboundary\_operator}}
        &   0.001 \textcolor{gray}{$\pm$ 0.087}
        &   0.008 \textcolor{gray}{$\pm$ 0.001}
        &   0.084 \textcolor{gray}{$\pm$ 0.006}
        &   \phantom{0}0.848 \textcolor{gray}{$\pm$ 0.011}
        &   \phantom{0}8.638 \textcolor{gray}{$\pm$ 0.045} \\
    \hline
        \hyperref[AlgNumModularVF]{\phantom{0}7. \textsf{num\_modular\_vf}}
        &   0.004 \textcolor{gray}{$\pm$ 0.754}
        &   0.030 \textcolor{gray}{$\pm$ 0.009}
        &   0.280 \textcolor{gray}{$\pm$ 0.001}
        &   \phantom{0}2.805 \textcolor{gray}{$\pm$ 0.016}
        &   \phantom{0}28.057 \textcolor{gray}{$\pm$ 0.107} \\
    \hline
        \hyperref[AlgNumCurl]{\phantom{0}8. \textsf{num\_curl\_operator}}
        &   0.022 \textcolor{gray}{$\pm$ 0.009}
        &   0.196 \textcolor{gray}{$\pm$ 0.024}
        &   1.923 \textcolor{gray}{$\pm$ 0.004}
        &   18.487 \textcolor{gray}{$\pm$ 0.136}
        &   182.774 \textcolor{gray}{$\pm$ 1.260} \\
    \hline
        \hyperref[AlgNumOneFormsB]{\phantom{0}9. \textsf{num\_one\_forms\_bracket}}
        &   0.058 \textcolor{gray}{$\pm$ 0.006}
        &   0.420 \textcolor{gray}{$\pm$ 0.007}
        &   4.278 \textcolor{gray}{$\pm$ 0.027}
        &   43.257 \textcolor{gray}{$\pm$ 0.071}
        &   434.450 \textcolor{gray}{$\pm$ 0.589} \\
    \hline
        \hyperref[AlgNumGauge]{10. \textsf{num\_gauge\_transformation}}
        &   0.051 \textcolor{gray}{$\pm$ 0.001}
        &   0.446 \textcolor{gray}{$\pm$ 0.010}
        &   4.380 \textcolor{gray}{$\pm$ 0.016}
        &   43.606 \textcolor{gray}{$\pm$ 0.212}
        &   434.704 \textcolor{gray}{$\pm$ 1.234} \\
    \hline
    \hline
 \end{tabular}
    }
\end{table}

To illustrate how fast the \textsf{NumPoissonGeometry} functions can be performed, we use the data in Table \ref{table:performance} to plot the time versus the number of points in each \ec{10^{\kappa}}--point (irregular) mesh on base $10$ \emph{log--log graphs}:
\begin{figure}[H]
    \centering
    \caption{Log-log graphs of the execution time in seconds versus the number of points in \ec{10^{\kappa}}--point (irregular) meshes of the \textsf{NumPoissonGeometry} functions 1--10 in Table \ref{table:performanceR2}, for \ec{\kappa=3,\ldots,7}. In red, the fitted linear model used to  predict the asymptotic behavior of the runtime for each function, with the corresponding determination coefficient (R-squared) indicated in each legend. We include a zoom-graph in each plot due to the accumulation of runtime values.} \label{fig:timeperformanceR2}
    \resizebox{\textwidth}{!}{
 \begin{tabular}{cccc}
        \includegraphics[width=\textwidth]{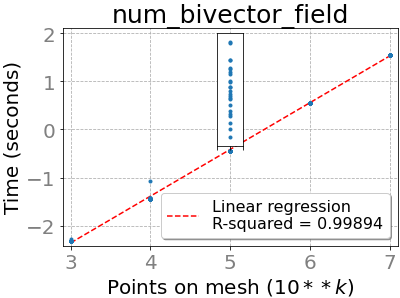} &
        \includegraphics[width=\textwidth]{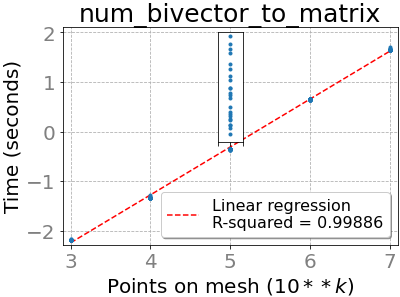} &
        \includegraphics[width=\textwidth]{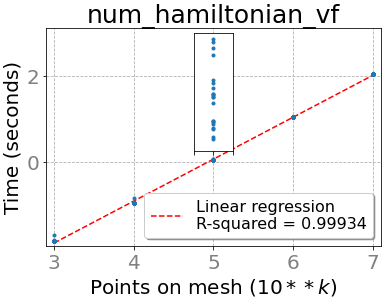} &
        \includegraphics[width=\textwidth]{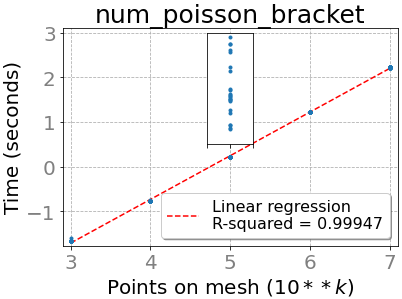} \\[0.5cm]
        \includegraphics[width=\textwidth]{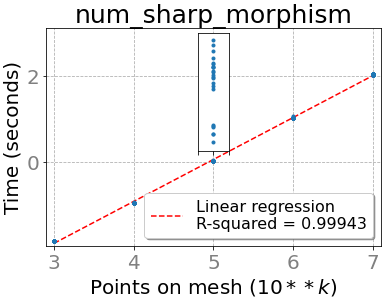} &
        \includegraphics[width=\textwidth]{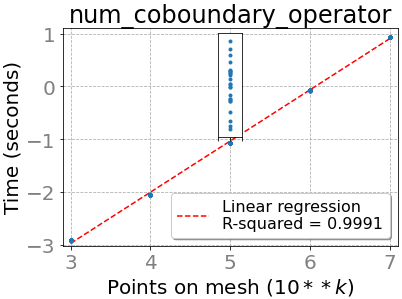} &
        \includegraphics[width=\textwidth]{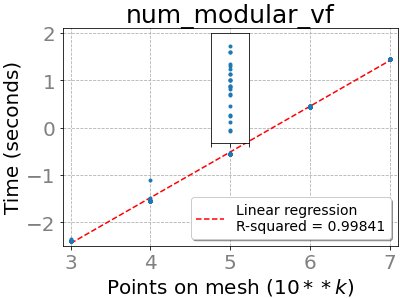} &
        \includegraphics[width=\textwidth]{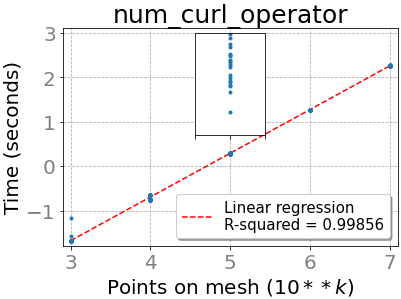} \\[0.5cm]
 \end{tabular}
    }
\end{figure}
\begin{figure}[H]
    \centering
    \resizebox{0.55\textwidth}{!}{
 \begin{tabular}{cccc}
        \multicolumn{2}{c}{\includegraphics[width=\textwidth]{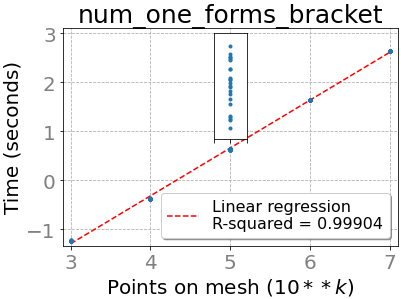}} &
        \multicolumn{2}{c}{\includegraphics[width=\textwidth]{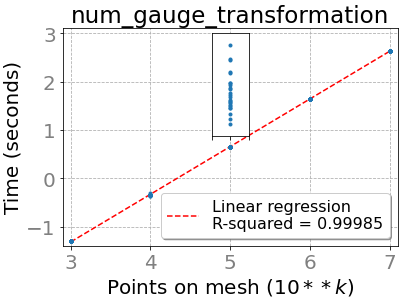}}
 \end{tabular}
    }
\end{figure}

We deduce from the graphs in Figure \ref{fig:timeperformanceR2} that, for the input data in Table \ref{table:InputR2}, all of our methods in Table \ref{table:performanceR2} were executed experimentally in polynomial time. Power-law relationships appear as straight lines in a log-log graph. Therefore the degree of the polynomial complexities are deduced by fitting a linear model and estimating its coefficient, which we carry out to 0.99 accuracy.

    \subsubsection{Three Dimensional Case}

For the performance tests of our methods in dimension three, we have set the Lie--Poisson bivector field \ec{\Pi_{\mathfrak{sl}(2)}} in (\ref{EcPiSL2}),
    \begin{equation*}
        \Pi_{\mathfrak{sl}(2)} =
        -x_3\frac{\partial}{\partial x_1}\wedge \frac{\partial}{\partial x_2} -
        x_2\frac{\partial}{\partial x_1}\wedge \frac{\partial}{\partial x_3} +
        x_1 \frac{\partial}{\partial x_2}\wedge \frac{\partial}{\partial x_3},
    \end{equation*}
associated to the 3--dimensional Lie algebra \ec{\mathfrak{sl}(2)}.

\begin{table}[H]
    \centering
    \caption{Input data used for the time performance tests of functions 1-11 in Table \ref{table:performance}.} \label{table:InputR2}
    \resizebox{\textwidth}{!}{
 \begin{tabular}{|c|l|}
        \hline
    \multicolumn{1}{|c|} {\textbf{Function}} & \multicolumn{1}{c|}{\textbf{Input}} \\
        \hline
        \hline
    1 & $\phantom{-}\Pi_{\mathfrak{sl}(2)}$ \\
         \hline
    2 & $\phantom{-}\Pi_{\mathfrak{sl}(2)}$ \\
         \hline
    3 & $\phantom{-}\Pi_{\mathfrak{sl}(2)}$, \ec{h = x_{1}^{2} + x_{2}^{2} - x_{3}^{2}} \\
         \hline
    4 & $\phantom{-}\Pi_{\mathfrak{sl}(2)}$, \ec{f = x_{1}^{2} + x_{2}^{2} - x_{3}^{2}}, \ec{g=x_{1} + x_{2} + x_{3}} \\
         \hline
    5 & $\phantom{-}\Pi_{\mathfrak{sl}(2)}$, \ec{\alpha = x_{1}\dd{x_{1}} + x_{2}\dd{x_{2}} - x_{3}\dd{x_{3}}} \\
         \hline
    6 & $\phantom{-}\Pi_{\mathfrak{sl}(2)}$, \ec{W = e^{{-1}/{(x_1^2 + x_2^2 - x_3^2)^2}} \big[ {x_1x_{3}}/(x_1^2 + x_2^2)\frac{\partial}{\partial{x_{1}}} + {x_2x_{3}}/(x_1^2 + x_2^2)\frac{\partial}{\partial{x_{2}}} + \frac{\partial}{\partial{x_{3}}} \big]} \\
         \hline
    7 & $\phantom{-}\Pi_{\mathfrak{sl}(2)}$, \ec{f=1} \\
         \hline
    8 & $\phantom{-}\Pi_{\mathfrak{sl}(2)}$, \ec{f=1} \\
         \hline
    9 & $\phantom{-}\Pi_{\mathfrak{sl}(2)}$, \ec{\alpha = x_{1}\dd{x_{1}} + x_{2}\dd{x_{2}} - x_{3}\dd{x_{3}}}, \ec{\beta = \dd{x_{1}} + \dd{x_{2}} + \dd{x_{3}}} \\
         \hline
    10 & $\phantom{-}\Pi_{\mathfrak{sl}(2)}$, \ec{\lambda = (x_{2}-x_{1})\dd{x_{1}} \wedge \dd{x_{2}} + (x_{3}-x_{1})\dd{x_{1}} \wedge \dd{x_{3}} + (x_{2}-x_{3})\dd{x_{2}} \wedge \dd{x_{3}}} \\
         \hline
    11 & $-\Pi_{\mathfrak{sl}(2)}$ \\
        \hline
        \hline
 \end{tabular}
    }
\end{table}

Table \ref{table:performance} lists the mean time in seconds (with standard deviation) it takes to evaluate the first eleven functions in \textsf{NumPoissonGeometry} (see, Table \ref{table:Funs-Algos-Exes}) on a irregular mesh on \ec{\RR{3}_{x}} with \ec{10^{\kappa}} points, computed by taking twenty-five samples, for \ec{\kappa = 3,\ldots,7}.

\begin{table}[H]
    \centering
    \caption{Summary of the time performance of \textsf{NumPoissonGeometry} functions in dimension 3.} \label{table:performance}
    \resizebox{\textwidth}{!}{
 \begin{tabular}{| l | l l l l l |}
    \hline
        \multirow{2}{*}{\hspace{1.25cm}\textbf{Function}} & \multicolumn{5}{c|}{\textbf{Points in mesh/Processing time (in seconds)}} \\
        & \hspace{0.85cm}$10^{3}$ & \hspace{0.85cm}$10^{4}$ & \hspace{0.85cm}$10^{5}$ & \hspace{0.85cm}$10^{6}$ & \hspace{0.85cm}$10^{7}$ \\
    \hline
    \hline
        \hyperref[AlgNumBivector]{\phantom{0}1. \textsf{num\_bivector\_field}}
        &   0.009 \textcolor{gray}{$\pm$ 0.009}
        &   0.051 \textcolor{gray}{$\pm$ 0.004}
        &   0.496 \textcolor{gray}{$\pm$ 0.002}
        &   \phantom{0}4.984 \textcolor{gray}{$\pm$ 0.023}
        &   \phantom{0}49.565 \textcolor{gray}{$\pm$ 0.222}  \\
    \hline
        \hyperref[AlgNumMatrixBivector]{\phantom{0}2. \textsf{num\_bivector\_to\_matrix}}
        &   0.008 \textcolor{gray}{$\pm$ 3.164}
        &   0.057 \textcolor{gray}{$\pm$ 0.002}
        &   0.553 \textcolor{gray}{$\pm$ 0.019}
        &   \phantom{0}5.442 \textcolor{gray}{$\pm$ 0.023}
        &   \phantom{0}55.249 \textcolor{gray}{$\pm$ 1.690} \\
    \hline
        \hyperref[AlgNumHamVF]{\phantom{0}3. \textsf{num\_hamiltonian\_vf}}
        &   0.017 \textcolor{gray}{$\pm$ 0.002}
        &   0.129 \textcolor{gray}{$\pm$ 0.001}
        &   1.263 \textcolor{gray}{$\pm$ 0.022}
        &   12.518 \textcolor{gray}{$\pm$ 0.064}
        &   126.091 \textcolor{gray}{$\pm$ 0.583} \\
    \hline
        \hyperref[AlgNumPoissonBracket]{\phantom{0}4. \textsf{num\_poisson\_bracket}}
        &   0.036 \textcolor{gray}{$\pm$ 0.001}
        &   0.299 \textcolor{gray}{$\pm$ 0.010}
        &   2.936 \textcolor{gray}{$\pm$ 0.067}
        &   29.600 \textcolor{gray}{$\pm$ 0.933}
        &   292.625 \textcolor{gray}{$\pm$ 6.094} \\
    \hline
        \hyperref[AlgNumSharp]{\phantom{0}5. \textsf{num\_sharp\_morphism}}
        &   0.017 \textcolor{gray}{$\pm$ 0.006}
        &   0.128 \textcolor{gray}{$\pm$ 0.005}
        &   1.252 \textcolor{gray}{$\pm$ 0.005}
        &   12.384 \textcolor{gray}{$\pm$ 0.038}
        &   124.851 \textcolor{gray}{$\pm$ 1.809} \\
    \hline
        \hyperref[AlgNumCoboundary]{\phantom{0}6. \textsf{num\_coboundary\_operator}}
        &   1.589 \textcolor{gray}{$\pm$ 0.016}
        &   1.705 \textcolor{gray}{$\pm$ 0.029}
        &   2.815 \textcolor{gray}{$\pm$ 0.032}
        &   12.972 \textcolor{gray}{$\pm$ 0.166}
        &   111.034 \textcolor{gray}{$\pm$ 1.365} \\
    \hline
        \hyperref[AlgNumModularVF]{\phantom{0}7. \textsf{num\_modular\_vf}}
        &   0.050 \textcolor{gray}{$\pm$ 0.001}
        &   0.103 \textcolor{gray}{$\pm$ 0.004}
        &   0.645 \textcolor{gray}{$\pm$ 0.006}
        &   \phantom{0}6.025 \textcolor{gray}{$\pm$ 0.013}
        &   \phantom{0}59.652 \textcolor{gray}{$\pm$ 0.146} \\
    \hline
        \hyperref[AlgNumCurl]{\phantom{0}8. \textsf{num\_curl\_operator}}
        &   0.019 \textcolor{gray}{$\pm$ 0.010}
        &   0.129 \textcolor{gray}{$\pm$ 0.027}
        &   1.199 \textcolor{gray}{$\pm$ 0.032}
        &   10.911 \textcolor{gray}{$\pm$ 0.181}
        &   105.841 \textcolor{gray}{$\pm$ 1.230} \\
    \hline
        \hyperref[AlgNumOneFormsB]{\phantom{0}9. \textsf{num\_one\_forms\_bracket}}
        &   0.093 \textcolor{gray}{$\pm$ 0.001}
        &   0.738 \textcolor{gray}{$\pm$ 0.007}
        &   7.285 \textcolor{gray}{$\pm$ 0.159}
        &   72.802 \textcolor{gray}{$\pm$ 1.474}
        &   724.514 \textcolor{gray}{$\pm$ 13.594} \\
    \hline
        \hyperref[AlgNumGauge]{10. \textsf{num\_gauge\_transformation}}
        &   0.051 \textcolor{gray}{$\pm$ 0.001}
        &   0.445 \textcolor{gray}{$\pm$ 0.010}
        &   4.395 \textcolor{gray}{$\pm$ 0.013}
        &   43.794 \textcolor{gray}{$\pm$ 0.173}
        &   437.326 \textcolor{gray}{$\pm$ 0.824} \\
    \hline
        \hyperref[AlgNumNormal]{11. \textsf{num\_linear\_normal\_form\_R3}}
        &   0.016 \textcolor{gray}{$\pm$ 0.438}
        &   0.061 \textcolor{gray}{$\pm$ 0.002}
        &   0.504 \textcolor{gray}{$\pm$ 0.012}
        &   \phantom{0}4.903 \textcolor{gray}{$\pm$ 0.017}
        &   \phantom{0}48.786 \textcolor{gray}{$\pm$ 0.219} \\
    \hline
    \hline
 \end{tabular}
    }
\end{table}

Now, to illustrate how fast the \textsf{NumPoissonGeometry} functions can be performed, we use the data in Table \ref{table:performance} to plot the time versus the number of points in each \ec{10^{\kappa}}--point (irregular) mesh on base $10$ \emph{log--log graphs}. These (execution time) plots are presented in Figures \ref{fig:timeperformance} and \ref{fig:timeperformance2}.

        \subsubsection{Polynomial Time}

\begin{figure}[H]
    \centering
    \caption{Log-log graphs of the execution time in seconds versus the number of points in \ec{10^{\kappa}}--point (irregular) meshes of the \textsf{NumPoissonGeometry} functions 1--5 and 8--11 in Table \ref{table:performance}, for \ec{\kappa=3,\ldots,7}. In red, the fitted linear model used to  predict the asymptotic behavior of the runtime for each function, with the corresponding determination coefficient (R-squared) indicated in each legend. We include a zoom-graph in each plot due to the accumulation of runtime values.}
    \label{fig:timeperformance}
    \resizebox{\textwidth}{!}{
 \begin{tabular}{cccc}
    \includegraphics[width=\textwidth]{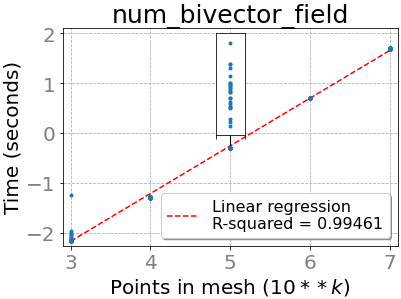} &
    \includegraphics[width=\textwidth]{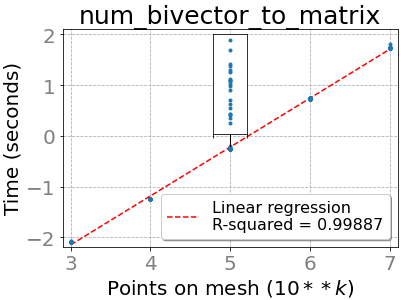} &
    \includegraphics[width=\textwidth]{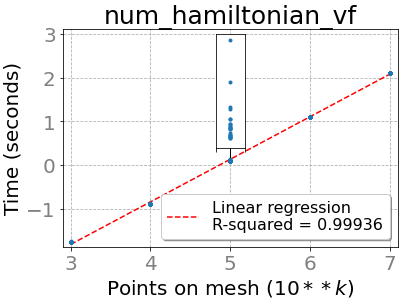} &
    \includegraphics[width=\textwidth]{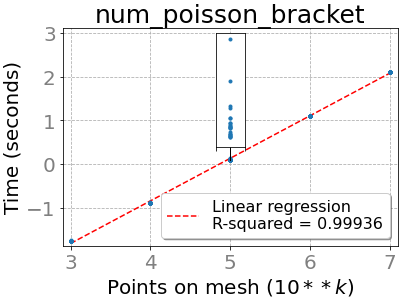} \\[0.5cm]
    \includegraphics[width=\textwidth]{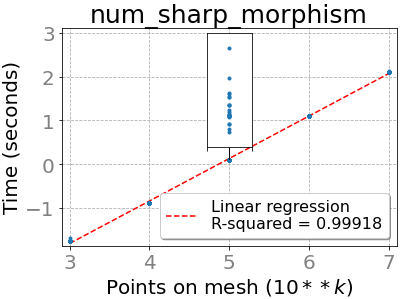} &
    \includegraphics[width=\textwidth]{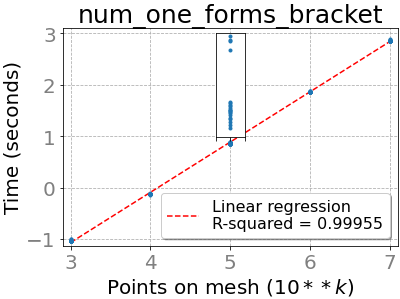} &
    \includegraphics[width=\textwidth]{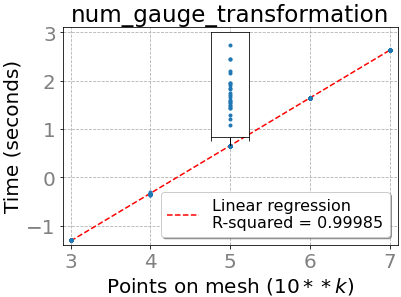} &
    \includegraphics[width=\textwidth]{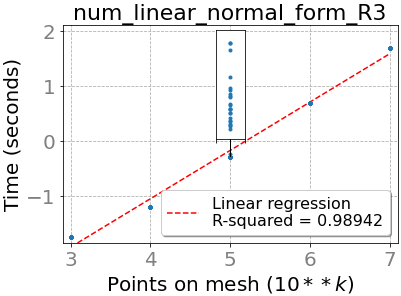} \\[0.5cm]
    \multicolumn{4}{c}{\includegraphics[width=\textwidth]{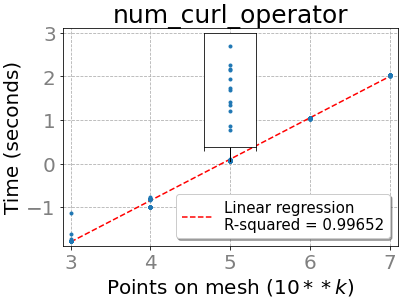}}
 \end{tabular}
    }
\end{figure}

The log-log graphs presented in Figure \ref{fig:timeperformance} correspond to \textsf{NumPoissonGeometry} methods that are executed in polynomial time of some degree. As power-law relationships appear as straight lines in a log-log graph, the complexities are deduced by fitting a linear model and estimating its coefficient to 0.99 accuracy.

Observe that the experimentally deduced time complexities of the \textsf{NumPoissonGeometry} functions in Figure \ref{fig:timeperformance} coincide with their theoretical time complexities in the following sense: the time complexities of methods 1--5 and 9--11 in Table \ref{table:complexity} are polynomially dependent on variables $m$ and $k$.

\begin{remark}
In our experiments, \textsf{num\_curl\_operator} approached a polynomial time complexity, although its theoretical time complexity is exponential (see, Table \ref{table:complexity}). This is not {\em a fortiori} a contradiction, because in Table \ref{table:complexity} we present an approximation of the \emph{worst--case} time complexity of the \textsf{NumPoissonGeometry} methods. In fact, it is an example that the execution time of our algorithms depends (naturally) on their inputs and that in some cases they run faster than expected.
\end{remark}

    \subsubsection{Exponential Time}

The log-log graphs presented in Figure \ref{fig:timeperformance} correspond to \textsf{NumPoissonGeometry} methods that are executed in exponential time. As exponential relationships trace polynomial curves in a log-log graph, the complexities are deduced by fitting a (non-linear) polynomial regression and estimating its coefficients to 0.99 accuracy.

\begin{figure}[H]
    \centering
    \caption{Log-log graphs of the execution time in seconds versus the number of points in \ec{10^{\kappa}}--point (irregular) meshes of the \textsf{NumPoissonGeometry} functions 6 and 7 in Table \ref{table:performance}, for \ec{\kappa=3,\ldots,7}. In red, the fitted model used to  predict the asymptotic behavior of the runtime for each function, with the corresponding determination coefficient (R-squared) indicated in each legend. We include a zoom-graph in each plot due to the accumulation of runtime values.}
    \label{fig:timeperformance2}
 \begin{tabular}{cc}
    \includegraphics[width=0.335\textwidth]{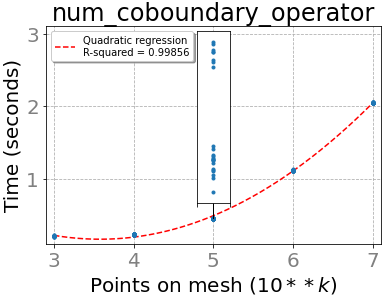} &
    \includegraphics[width=0.35\textwidth]{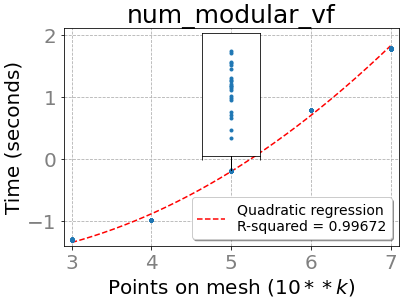}
 \end{tabular}
\end{figure}

Observe that although the \textsf{NumPoissonGeometry} functions presented in Figure \ref{fig:timeperformance2} are executed in exponential time, they are relatively fast since their evaluation on a 10 million point (irregular) mesh takes at most 2 minutes in the experiments on our desktop workstation. Furthermore, the execution times of the \textsf{num\_modular\_vf} function are fitted to a linear model in the interval [4,7], as we illustrate in the following figure:

\begin{figure}[H]
  \centering
  \caption{Log-log graph of the execution time in seconds versus the number of points in a \ec{10^{\kappa}}--point (irregular) mesh of the \textsf{NumPoissonGeometry} function \textsf{num\_modular\_vf}, for \ec{\kappa=4,\ldots,7}. In red, the fitted linear model used to predict the asymptotic behavior of the runtime, with the corresponding determination coefficient (R-squared) indicated in the legend.}
  \includegraphics[width=0.35\textwidth]{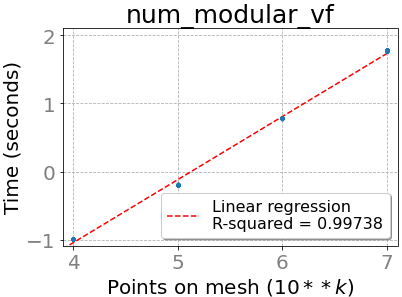}
\end{figure}

    \subsubsection{Flaschka--Ratiu Bivector Fields}

For the performance tests of the method \textsf{num\_flaschka\_ratiu\_bivector}, we have used as inputs the following scalar functions on \ec{\RR{4}_{x}} (see Example \ref{example:so3FR}):
    \begin{equation*}
        K_{1} = \tfrac{1}{2}x_{4}, \qquad K_{2} = -x_{1}^{2} + x_{2}^{2} + x_{3}^{2}
    \end{equation*}

\begin{table}[H]
    \centering
    \caption{Mean time in seconds (with standard deviation) it takes to evaluate the \textsf{num\_flaschka\_ratiu\_bivector} method on a irregular mesh on $\RR{4}$ with $10^{\kappa}$ points, computed by taking twenty-five samples, for \ec{\kappa = 3,\ldots,7}.} \label{table:performanceFR}
    \resizebox{\textwidth}{!}{
 \begin{tabular}{|l|l l l l l|}
    \hline

        \multirow{2}{*}{\hspace{1.25cm}\textbf{Function}} & \multicolumn{5}{c|}{\textbf{Points in mesh/Processing time (in seconds)}} \\
        & \hspace{0.85cm}$10^{3}$ & \hspace{0.85cm}$10^{4}$ & \hspace{0.85cm}$10^{5}$ & \hspace{0.85cm}$10^{6}$ & \hspace{0.85cm}$10^{7}$ \\
    \hline
    \hline
        \hyperref[AlgNumFRatiu]{\textsf{num\_flaschka\_ratiu\_bivector}}
        &   0.0158 \textcolor{gray}{$\pm$ 0.105}
        &   0.057 \textcolor{gray}{$\pm$  0.003}
        &   0.505 \textcolor{gray}{$\pm$ 0.003}
        &   \phantom{0}4.993 \textcolor{gray}{$\pm$ 0.029}
        &   \phantom{0}49.563 \textcolor{gray}{$\pm$ 0.207}  \\
    \hline
    \hline
 \end{tabular}
    }
\end{table}

Figure \ref{fig:timeperformanceFR} below illustrates that the \textsf{num\_flaschka\_ratiu\_bivector} method is executed in polynomial time of some degree, which coincides with its theoretical complexity presented in Table \ref{table:complexity} in the sense that it is polynomially dependent on variables $m$ and $k$. As power-law relationships appear as straight lines in a log-log graph, the complexity is deduced by fitting a linear model and estimating its coefficient to 0.99 accuracy.

\begin{figure}[H]
  \centering
    \caption{Log-log graph of the execution time in seconds versus the number of points in a \ec{10^{\kappa}}--point (irregular) mesh of the \textsf{num\_flaschka\_ratiu\_bivector} method, for \ec{\kappa=3,\ldots,7}. In red, the fitted linear model used to predict the asymptotic behavior of the runtime, with the corresponding determination coefficient (R-squared) indicated in the legend. We include a zoom-graph due to the accumulation of runtime values.}
  \label{fig:timeperformanceFR}
  \includegraphics[width=0.35\textwidth]{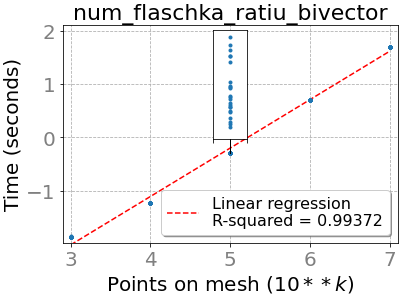}
\end{figure}

    \section{Appendix}

In this section we explain how mesheses can be used to evaluate our methods in \textsf{NumPoissonGeometry}. Our module does not generate meshes by default, this is to prevent the complexity of our algorithms from increasing and the overall computational performance from decreasing.

A mesh can be easily written as a \textsf{NumPy} array (see, Subsection \ref{subsec:mesh}). Recall that \textsf{NumPy} arrays can be generated with evenly distributed points, or with randomly distributed points as in the main text above and also shown in the next examples below. They could also come from preloaded data points of interest.

The mesh \ec{Q = \{[1,0,0], [0,1,0], [0,0,1]\} \subset \RR{3}} can be written as a \textsf{NumPy} array as follows:
\begin{tcolorbox}[arc=0mm, boxsep=0mm, skin=bicolor, colback=pink!15, colframe=blue!20, colbacklower=blue!0, breakable, halign=left]
\LineCode \textsf{points = [[1, 0, 0], [0, 1, 0], [0, 0, 1]]}
    \hspace*{\fill} \CommentCode{a list of lists encoding $Q$} \\
\LineCode \textsf{mesh = numpy.array(points)} \hspace*{\fill} \CommentCode{convert \textsf{points} into a \textsf{NumPy} array}
\end{tcolorbox}

To generate a random mesh with \ec{10^{6}} points in \ec{\RR{3}} as an \textsf{NumPy} array, we can run the following:
\begin{tcolorbox}[arc=0mm, boxsep=0mm, skin=bicolor, colback=pink!15, colframe=blue!20, colbacklower=blue!0, breakable, halign=left]
\LineCode \textsf{mesh = numpy.random.rand(10**6, 3)} \hspace*{\fill} \CommentCode{\ec{(10^{6}, 3)} \textsf{NumPy} array with random samples from a uniform distribution over [0,1)}
\end{tcolorbox}

    \subsection*{Acknowledgements}

This research was partially supported by CONACyT and UNAM-DGAPA-PAPIIT-IN104819. JCRP thanks CONACyT for a postdoctoral fellowship held during the production of this work.


\end{document}